\documentclass[leqno,11pt]{amsart}
\usepackage{amsmath,amscd,amsthm}
\usepackage{epsfig,graphics,ifsym,color}
\usepackage{amssymb,latexsym}
\usepackage{mathrsfs,eucal}
\usepackage{hyperref}
\usepackage[all,cmtip]{xy}
\usepackage{arydshln}

\newtheorem{thm}{\bf Theorem}[section]
\newtheorem{df}[thm]{\bf Definition}
\newtheorem{prop}[thm]{\bf Proposition}
\newtheorem{cor}[thm]{\bf Corollary}
\newtheorem{lem}[thm]{\bf Lemma}
\newtheorem{rem}[thm]{\bf Remark}
\newtheorem{ex}[thm]{\bf Example}

\newtheorem{nono-theorem}{Theorem}[]
\newtheorem*{thm*}{Theorem}

\newcommand{\B}{\mathbf{B}}

\newcommand{\cP}{\mathscr{P}}

\newcommand{\pf}{\noindent{\bfseries Proof. }}
\newcommand{\ov}{\overline}

\newcommand{\hf}{\frac{1}{2}}
\newcommand{\N}{\mathbb{N}}
\newcommand{\gl}{\mathfrak{gl}}
\newcommand{\Z}{\mathbb{Z}}

\newcommand{\te}{\widetilde{e}}
\newcommand{\tf}{\widetilde{f}}
\newcommand{\g}{\mathfrak{g}}
\newcommand{\td}{\widetilde}

\newcommand{\mc}{\mathcal}
\newcommand{\mf}{\mathfrak}
\newcommand{\La}{\Lambda}

\newcommand{\la}{\lambda}

\numberwithin{equation}{section}

\begin{document}
\title[Combinatorial extension of stable branching rules]
{Combinatorial extension of stable branching rules for classical groups}
\author{JAE-HOON KWON}

\address{Department of Mathematical Sciences, Seoul National University, Seoul 08826, Korea}
\email{jaehoonkw@snu.ac.kr}

\keywords{branching rules, classical groups, quantum groups, crystal graphs}
\subjclass[2010]{17B37, 22E46, 05E10}

\thanks{This work was supported by Samsung Science and Technology Foundation under Project Number SSTF-BA1501-01.}

\begin{abstract}
We give new combinatorial formulas for decomposition of the tensor product of integrable highest weight modules over the classical Lie algebras of type $B, C, D$, and the branching decomposition of an integrable highest weight module with respect to a maximal Levi subalgebra of type $A$. This formula is based on a combinatorial model of classical crystals called spinor model.
We show that our formulas extend in a bijective way various stable branching rules for classical groups to arbitrary highest weights, including the Littlewood restriction rules. 
\end{abstract}

\maketitle
\setcounter{tocdepth}{1}

\section{Introduction}

\subsection{Stable branching rules}
Let $G$ be a classical group over the complex numbers with  a closed subgroup $H$.
For a finite-dimensional irreducible $G$-module $V$, there are various branching rules for the multiplicity of an irreducible $H$-module in $V$, which are given as a sum of product of Littlewood-Richardson coefficients. 
These formulas, which originated in the Littlewood restriction rules \cite{Lw-1,Lw-2}, are often referred to as {\em stable branching rules} since they hold  in a certain range of highest weights of irreducible modules and depend only on the partitions parameterizing highest weights when the ranks of classical groups are sufficiently large. We refer the reader to \cite{HTW} for a systematic approach to these formulas and detailed exposition on previous works.

The goal of this paper is to give combinatorial extension of various stable branching rules to arbitrary highest weights, that is, to give a combinatorial formula for the branching multiplicity, which holds for arbitrary highest weights and also extends a given stable branching rule in a bijective way.

We recall that there are formulas extending stable branching rules to arbitrary highest weights (for example, \cite{EW,Kng-b,Ko,KoT87,KoT90,Ne,We}) where the multiplicities outside a stable range are in general given as an alternating sum of products of Littlewood-Richardson coefficients or determined by certain modification and cancellation of irreducible factors in the decomposition inside a stable range. But compared to these results which are obtained in an algebraic way, relatively not much is known about combinatorial or subtraction-free extension. 
On the other hand, the theory of crystals \cite{Kas1,Kas94} has made it possible to provide formulas for the branching multiplicities, which are available for arbitrary highest weights and described explicitly in terms of various combinatorial objects (for example, \cite{BZ-2,Kam,KashNaka,Li}). However, most of the known results do not seem to imply the stable branching rules immediately. So one may ask which combinatorial model for classical crystals explains and extends the stable branching formulas more naturally.

\subsection{Combinatorial extensions}\label{subsec:main result}

Let us explain our results in details. Our approach is based on the theory of crystal base, and the theory of reductive dual pairs \cite{H}.

Let $\mf g_\infty$ be the affine Kac-Moody Lie algebra of type $B_\infty$, $C_\infty$, and $D_\infty$ \cite{Kac90}, and let $\mf {l}_\infty$ be its maximal Levi subalgebra of type $A$.
We first
give a new formula for decomposition of the tensor product of integrable highest weight $\mf g_\infty$-modules, and the branching decomposition of an integrable highest weight $\mf g_\infty$-module with respect to  $\mf {l}_\infty$. 

Our formula is given in terms of a combinatorial model for the $\mf g_\infty$-crystal associated to an integrable highest weight $\mf g_\infty$-module, which was recently introduced by the author \cite{K13,K14}.
A main advantage of this model, which we call {\em spinor model}, is that it is compatible with the combinatorics of usual Young and Littlewood-Richardson tableaux, and hence it turns out to fit into the stability phenomenon of decomposition numbers very nicely.

The notion of spinor model is motivated by the duality on the fermionic Fock space $\mc{F}^{\frac{n}{2}}$, where $(\mf{g}_\infty, {\rm G}_n)$ forms a reductive dual pair on $\mc{F}^{\frac{n}{2}}$ for some complex classical reductive algebraic group ${\rm G}_n$ \cite{H,Wa}. Let $\mc{P}({\rm G}_n)$ be the set of partitions parameterizing the highest weights for the finite-dimensional irreducible ${\rm G}_n$-modules appearing in $\mc{F}^{\frac{n}{2}}$. By dual correspondence, we may parameterize the dominant integral weights for $\mf{g}_\infty$ by  $\mc{P}({\rm G}_n)$ for all $n$. Let ${\bf T}^{\mf g}(\la,n)$ denote the spinor model of the crystal associated to an integrable highest weight $\mf{g}_\infty$-module with highest weight corresponding to $\la\in \mc{P}({\rm G}_n)$. Here $\mf{g}=\mf{b}, \mf{c}, \mf{d}$ is understood as a symbol representing the type of $\mf g_\infty$.

Now, let $\la\in \mc{P}({\rm G}_{m+n})$, $\mu\in {\mc P}({\rm G}_m)$, and $\nu\in {\mc P}({\rm G}_{n})$ be given. The multiplicity of ${\bf T}^{\mf g}(\la,m+n)$ in the tensor product ${\bf T}^{\mf g}(\mu,m)\otimes {\bf T}^{\mf g}(\nu,n)$ is given by the cardinality of a subset ${\tt \large LR}^{\la}_{\mu\nu }({\mf g})$ of ${\bf T}^{\mf g}(\nu,n)$ satisfying certain combinatorial conditions determined by the general theory of crystal base (Proposition \ref{prop:tensor product}). Hence we may regard elements in ${\tt \large LR}^{\la}_{\mu\nu }({\mf g})$ as a generalization of Littlewood-Richardson tableaux. Then we show as our first main result that ${\tt \large LR}^{\la}_{\mu\nu }({\mf g})$ has a nice decomposition into a set of pairs of Littlewood-Richardson tableaux in a stable range. More precisely, we construct an explicit bijection (Theorem \ref{thm:stable tensor})
\begin{equation}\label{eq:1st main result}
{\tt \large LR}^{\la}_{\mu\nu }({\mf g}) \longrightarrow 
\bigsqcup_{\gamma\in\cP}\bigsqcup_{\delta\in \cP_{\mf g}} {\tt LR}^{\gamma}_{\mu' \delta}\times {\tt LR}^{\la'}_{\gamma\nu'}\, ,
\quad\quad\text{when $\ell(\la)\leq \frac{1}{2}\min\{m,n\}$}.
\end{equation}
Here $\cP$ is the set of partitions, $\cP_{\mf g}$ is a subset of $\cP$ given in \eqref{partition of type x}, ${\tt LR}^{\alpha}_{\beta \gamma}$ is the set of Littlewood-Richardson tableaux of shape $\alpha/\beta$ and content $\gamma$ for $\alpha, \beta, \gamma\in\cP$, $\la'$ is the conjugate of $\la$, and $\ell(\la)$ is the length of $\la$.

Next, we consider the branching decomposition with respect to $\mf l_\infty$.
Let $\mu\in {\mc P}({\rm G}_n)$ be given. For $\la\in \cP$, we denote by $\texttt{\large LR}^{\mu}_{\la}({\mf g})$ the subset of ${\bf T}^{\mf g}(\mu,n)$, which consists of the highest weight elements of the $\mf l_\infty$-subcrystals in ${\bf T}^{\mf g}(\mu,n)$ isomorphic to the crystal of Young tableaux of shape $\la'$. The cardinality of $\texttt{\large LR}^{\mu}_{\la}({\mf g})$ is the associated branching multiplicity.
We show as our second main result that $\texttt{\large LR}^{\mu}_{\la}({\mf g})$ decomposes into a set of Littlewood-Richardson tableaux in a stable range by constructing an explicit bijection
\begin{equation}\label{eq:2nd main result}
{\tt \large LR}^\mu_{\la}({\mf g}) \longrightarrow  \bigsqcup_{\delta\in \cP_{\mf g}} {\tt LR}^{\la'}_{\delta\mu'}
\quad\quad\text{when $\ell(\la) \leq \frac{n}{2}$}.
\end{equation}

Finally, by the reciprocity laws associated to suitable see-saw pairs on $\mc{F}^{\frac{n}{2}}$, the bijections in \eqref{eq:1st main result} and \eqref{eq:2nd main result} recover well-known stable branching rule for the pairs

\begin{itemize}
\item[(I)] ${\rm G}_{m+n}\supset {\rm G}_m\times {\rm G}_n$,\quad 
\item[(II)]  ${\rm GL}_n \supset {\rm G}_n$   \ (the Littlewood restriction rule),

\end{itemize}
with ${\rm G}_\ell={\rm Sp}_\ell$, ${\rm O}_\ell$ (Theorems \ref{thm:stable tensor-dual} and \ref{thm:stable branching-dual}), and therefore the sets ${\tt \large LR}^{\la}_{\mu\nu }({\mf g})$ and ${\tt \large LR}^\mu_{\la}({\mf g})$ extend these stable branching rules to arbitrary highest weights in a bijective way.
On the other hand, considering a subset of tableaux in ${\bf T}^{\mf g}(\la)$ with entries from $\{1,\ldots,k\}$, we can also describe the tensor product decomposition and branching decomposition for crystals of type $B_k$, $C_k$, and $D_k$ in terms of ${\tt \large LR}^{\la}_{\mu\nu }({\mf g})$ and ${\tt \large LR}^\mu_{\la}({\mf g})$. That is, we have branching rules for the pairs
\begin{itemize}
\item[(III)]  ${\rm G}_\ell \times {\rm G}_\ell  \supset {\rm G}_\ell$,\quad  \item[(IV)] ${\rm G}_\ell \supset {\rm GL}_{[\ell/2]}$,
\end{itemize}
with ${\rm G}_\ell={\rm Sp}_\ell$, ${\rm Spin}_\ell$, ${\rm Osp}_{1|\ell}$ and $k=[\ell/2]$, together with the stable limits \eqref{eq:1st main result} and \eqref{eq:2nd main result}, respectively (see Remark \ref{rem:reduction to k}). In this case, the formula \eqref{eq:1st main result} for (III) looks different from the known result, and the formulas in (III) and (IV) for ${\rm Spin}_\ell$-modules with half-integer weights seem to be new.

An extension of stable branching formula similar to \eqref{eq:1st main result} and \eqref{eq:2nd main result} has been obtained in \cite{K12} by using a different model of $\mf{g}_\infty$-crystals for $\mf g=\mf b$, $\mf c$ but not for $\mf g=\mf d$, which is a technically more difficult case. For other combinatorial extensions, we refer the reader to \cite{Su} for (II) with ${\rm G}_n={\rm Sp}_n$, and \cite[Proposition 2.6]{BZ-1} for (IV) with ${\rm G}_\ell={\rm Sp}_\ell$, ${\rm SO}_\ell$.

Recently, there have been several works studying stability phenomenon from completely different  viewpoints from ours.
In \cite{DPS,PS,SS} (see also references therein), the branching multiplicities in a stable range are studied in the context of a non-semisimple tensor category of modules over classical Lie algebras of infinite rank. Also in \cite{We}, the extension of stable branching multiplicities in (II) has been studied using $q$-versions of Brauer’s centralizer algebras.

We expect that our combinatorial approach would give a new insight into the structure of representations of classical Lie algebras and groups which are related to the stability phenomenon.


\subsection{Application to holomorphic discrete series}
We give an interesting application of spinor model to another stability phenomenon.

In \cite{HTW} a systematic approach of stable branching rule is given by using Howe's theory of reductive dual pairs \cite{H}, which also partly motivated the work in the present paper.
Instead of integrable highest weight $\mf g_\infty$-modules, a family of infinite-dimensional unitarizable representations of classical Lie algebras of finite rank, which appears in the study of dual pair acting on bosonic spaces, is considered in \cite{HTW}, and then the irreducibility of their associated generalized Verma modules in a certain stable range plays a crucial role in giving a unified proof of 10 families of stable branching formulas classified there.

We show that the spinor model ${\bf T}^\g(\la,n)$ admits a variation to give a combinatorial character formula for these unitarizable representations (Theorem \ref{thm:char of unitary}) thanks to super duality \cite{CLW}. Hence we obtain the tensor product multiplicity and the branching multiplicity with respect to a maximal Levi subalgebra type $A$ for these irreducible modules. In particular, the branching multiplicity formula implies the irreducibility of associated generalized Verma modules in a stable range called {\em holomorphic discrete series} (Corollary \ref{cor:stable branching for unitary}).

\subsection{Organization}
The paper is organized as follows. In Section 2, we briefly recall necessary background. In Section 3, we review the spinor model of crystals of types $BCD$. In Section 4, we establish the bijections \eqref{eq:1st main result} and \eqref{eq:2nd main result}. In Section 5, we recover well-known stable branching rules by using the dualities on Fock spaces. In Section 6, we discuss a combinatorial character formula for infinite-dimensional unitarizable representations.

\section{Preliminaries}

\subsection{Lie algebras of type $B$, $C$, and $D$}\label{subsec:notations}
We assume that our base field is $\mathbb{C}$, and $\Z_+$ (resp. $\N$) denotes the set of non-negative (resp. positive) integers.
Let  $\mf{b}_{\infty}$, $\mf{c}_{\infty}$, and $\mf{d}_{\infty}$ be  the Kac-Moody Lie algebras of type $B_\infty$, $C_\infty$, and $D_\infty$, respectively (see \cite{Kac90}). Following  the conventions of $\mf{g}_\infty$ (${\mf g}={\mf b}, {\mf c}, {\mf d}$) in \cite[Section 2.2]{CK}, we let
\begin{itemize}

\item[$\cdot$] $I=\Z_+$ : the index set for simple roots,

\item[$\cdot$] $\Pi=\{\,\alpha_i\,|\,i\in I\,\}$ : the set of
simple roots,

\item[$\cdot$] $\Pi^{\vee}=\{\,\alpha_i^{\vee}\,|\,i\in I\,\}$ : the set of
simple coroots,

\item[$\cdot$] $\Lambda_i^{\mf{g}}$ : the fundamental weight with respect to $\alpha_i$,

\item[$\cdot$] $P=\Z\Lambda_0^{\mf g}\oplus \bigoplus_{i\geq 1}\Z\epsilon_i$ : the weight lattice,
\end{itemize}
where the Dynkin diagram of the Cartan matrix $(\langle \alpha_j,\alpha^\vee_i \rangle)$, $\Pi$, and  $\La^{\mf g}_i$ are given by
\begin{center}
\hskip -1cm  \setlength{\unitlength}{0.16in}
\begin{picture}(24,4)
\put(2,2){\makebox(0,0)[c]{${\mf b}_{\infty}\ :$}}

\put(5.6,2){\makebox(0,0)[c]{$\bigcirc$}}
\put(8,2){\makebox(0,0)[c]{$\bigcirc$}}
\put(10.4,2){\makebox(0,0)[c]{$\bigcirc$}}
\put(14.85,2){\makebox(0,0)[c]{$\bigcirc$}}
\put(17.25,2){\makebox(0,0)[c]{$\bigcirc$}}
\put(19.4,2){\makebox(0,0)[c]{$\bigcirc$}}
\put(8.35,2){\line(1,0){1.5}} \put(10.82,2){\line(1,0){0.8}}
\put(13.2,2){\line(1,0){1.2}} \put(15.28,2){\line(1,0){1.45}}
\put(17.7,2){\line(1,0){1.25}} \put(19.81,2){\line(1,0){0.9}}
\put(6.8,1.97){\makebox(0,0)[c]{$\Longleftarrow$}}
\put(12.5,1.95){\makebox(0,0)[c]{$\cdots$}}
\put(21.5,1.95){\makebox(0,0)[c]{$\cdots$}}
\put(5.6,1){\makebox(0,0)[c]{\tiny ${\alpha}_0$}}
\put(8,1){\makebox(0,0)[c]{\tiny ${\alpha}_1$}}
\put(10.4,1){\makebox(0,0)[c]{\tiny ${\alpha}_2$}}
\put(15,1){\makebox(0,0)[c]{\tiny ${\alpha}_{k-1}$}}
\put(17.4,1){\makebox(0,0)[c]{\tiny ${\alpha}_k$}}
\put(19.8,1){\makebox(0,0)[c]{\tiny ${\alpha}_{k+1}$}}
\end{picture}
\begin{equation*}
\begin{split}
&\Pi=\{\, \alpha_0=-\epsilon_1, \
 \alpha_i=\epsilon_i-\epsilon_{i+1} \ (i\geq 1)\, \},\\
&\La^{\mf b}_i =2\La^{\mf b}_0+\epsilon_1+\cdots+\epsilon_i, \ \ \
(i\geq 1),
\end{split}
\end{equation*}
\end{center}

\begin{center}
\hskip -1cm  \setlength{\unitlength}{0.16in}
\begin{picture}(24,4)
\put(2,2){\makebox(0,0)[c]{${\mf c}_{\infty}\ :$}}

\put(5.6,2){\makebox(0,0)[c]{$\bigcirc$}}
\put(8,2){\makebox(0,0)[c]{$\bigcirc$}}
\put(10.4,2){\makebox(0,0)[c]{$\bigcirc$}}
\put(14.85,2){\makebox(0,0)[c]{$\bigcirc$}}
\put(17.25,2){\makebox(0,0)[c]{$\bigcirc$}}
\put(19.4,2){\makebox(0,0)[c]{$\bigcirc$}}
\put(8.35,2){\line(1,0){1.5}} \put(10.82,2){\line(1,0){0.8}}
\put(13.2,2){\line(1,0){1.2}} \put(15.28,2){\line(1,0){1.45}}
\put(17.7,2){\line(1,0){1.25}} \put(19.81,2){\line(1,0){0.9}}
\put(6.8,1.97){\makebox(0,0)[c]{$\Longrightarrow$}}
\put(12.5,1.95){\makebox(0,0)[c]{$\cdots$}}
\put(21.5,1.95){\makebox(0,0)[c]{$\cdots$}}
\put(5.6,1){\makebox(0,0)[c]{\tiny ${\alpha}_0$}}
\put(8,1){\makebox(0,0)[c]{\tiny ${\alpha}_1$}}
\put(10.4,1){\makebox(0,0)[c]{\tiny ${\alpha}_2$}}
\put(15,1){\makebox(0,0)[c]{\tiny ${\alpha}_{k-1}$}}
\put(17.4,1){\makebox(0,0)[c]{\tiny ${\alpha}_k$}}
\put(19.8,1){\makebox(0,0)[c]{\tiny ${\alpha}_{k+1}$}}

\end{picture}
\begin{equation*}
\begin{split}
&\Pi=\{\,\alpha_0=-2\epsilon_1 , \ \alpha_i=\epsilon_i-\epsilon_{i+1} \ (i\geq 1)\, \},\\
&\La^{\mf c}_i=\La^{\mf c}_0+\epsilon_1+\cdots+\epsilon_i, \ \ \ \ \
(i\geq 1),
\end{split}
\end{equation*}
\end{center}\vskip 5mm

\begin{center}
\setlength{\unitlength}{0.16in} \hskip -3cm
\hskip 2cm \begin{picture}(24,5.8)
\put(2,2){\makebox(0,0)[c]{${\mf d}_{\infty}\ :$}}
\put(6,0){\makebox(0,0)[c]{$\bigcirc$}}
\put(6,4){\makebox(0,0)[c]{$\bigcirc$}}
\put(8,2){\makebox(0,0)[c]{$\bigcirc$}}
\put(10.4,2){\makebox(0,0)[c]{$\bigcirc$}}
\put(14.85,2){\makebox(0,0)[c]{$\bigcirc$}}
\put(17.25,2){\makebox(0,0)[c]{$\bigcirc$}}
\put(19.4,2){\makebox(0,0)[c]{$\bigcirc$}}
%
\put(6.35,0.3){\line(1,1){1.35}} \put(6.35,3.7){\line(1,-1){1.35}}
\put(8.4,2){\line(1,0){1.55}} \put(10.82,2){\line(1,0){0.8}}
\put(13.2,2){\line(1,0){1.2}} \put(15.28,2){\line(1,0){1.45}}
\put(17.7,2){\line(1,0){1.25}} \put(19.8,2){\line(1,0){1.25}}
%
\put(12.5,1.95){\makebox(0,0)[c]{$\cdots$}}
\put(22,1.95){\makebox(0,0)[c]{$\cdots$}}
\put(6,5){\makebox(0,0)[c]{\tiny $\alpha_{0}$}}
\put(6,-1.2){\makebox(0,0)[c]{\tiny $\alpha_{1}$}}
\put(8.2,1){\makebox(0,0)[c]{\tiny $\alpha_{2}$}}
\put(10.4,1){\makebox(0,0)[c]{\tiny $\alpha_{3}$}}
\put(14.8,1){\makebox(0,0)[c]{\tiny $\alpha_{k-1}$}}
\put(17.15,1){\makebox(0,0)[c]{\tiny $\alpha_k$}}
\put(19.5,1){\makebox(0,0)[c]{\tiny $\alpha_{k+1}$}}
\end{picture}

\begin{equation*}
\begin{split}
&\Pi=\{\, \alpha_0=-\epsilon_1-\epsilon_2, \  \alpha_i=\epsilon_i-\epsilon_{i+1} \ (i\geq 1)\, \},\\
&\La^{\mf d}_i =
\begin{cases}
\La^{\mf d}_0+\epsilon_1, & \text{if $i=1$}, \\
2\La^{\mf d}_0+\epsilon_1+\cdots+\epsilon_i, & \text{if $i>1$}.
\end{cases}
\end{split}
\end{equation*}
\end{center}\vskip 5mm
We also let $\mathfrak{l}_\infty$ be the subalgebra of ${\mf g}_\infty$ associated to $\{\,\alpha_i\,|\, i\in I\setminus \{0\}\,\}$,  which is of type $A_{+\infty}$ (cf. \cite{Kac90}). We assume that its weight lattice is $\bigoplus_{i\geq 1}\Z\epsilon_i \subset P$.

We will also consider the following Kac-Moody Lie superalgebra $\mf{b}^{\bullet}_\infty$ of infinite rank, whose Dynkin diagram, $\Pi$, and the fundamental weights $\La^{\mf b^\bullet}_i$ are given by

\begin{center}
\hskip -1cm  \setlength{\unitlength}{0.16in}
\begin{picture}(24,4)
\put(2,2){\makebox(0,0)[c]{${\mf b}^\bullet_{\infty}\ :$}}

\put(5.6,2){\circle*{0.9}}
\put(8,2){\makebox(0,0)[c]{$\bigcirc$}}
\put(10.4,2){\makebox(0,0)[c]{$\bigcirc$}}
\put(14.85,2){\makebox(0,0)[c]{$\bigcirc$}}
\put(17.25,2){\makebox(0,0)[c]{$\bigcirc$}}
\put(19.4,2){\makebox(0,0)[c]{$\bigcirc$}}
\put(8.35,2){\line(1,0){1.5}} \put(10.82,2){\line(1,0){0.8}}
\put(13.2,2){\line(1,0){1.2}} \put(15.28,2){\line(1,0){1.45}}
\put(17.7,2){\line(1,0){1.25}} \put(19.81,2){\line(1,0){0.9}}
\put(6.8,1.97){\makebox(0,0)[c]{$\Longleftarrow$}}
\put(12.5,1.95){\makebox(0,0)[c]{$\cdots$}}
\put(21.5,1.95){\makebox(0,0)[c]{$\cdots$}}
\put(5.6,1){\makebox(0,0)[c]{\tiny ${\alpha}_0$}}
\put(8,1){\makebox(0,0)[c]{\tiny ${\alpha}_1$}}
\put(10.4,1){\makebox(0,0)[c]{\tiny ${\alpha}_2$}}
\put(15,1){\makebox(0,0)[c]{\tiny ${\alpha}_{k-1}$}}
\put(17.4,1){\makebox(0,0)[c]{\tiny ${\alpha}_k$}}
\put(19.8,1){\makebox(0,0)[c]{\tiny ${\alpha}_{k+1}$}}
\end{picture}
\begin{equation*}
\begin{split}
&\Pi=\{\, \alpha_0=-\epsilon_1, \
 \alpha_i=\epsilon_i-\epsilon_{i+1} \ (i\geq 1)\, \},\\
&\La^{\mf b^\bullet}_i =2\La^{\mf b}_0+\epsilon_1+\cdots+\epsilon_i, \ \ \
(i\geq 1),
\end{split}
\end{equation*}
\end{center}
where\ \  $
\begin{picture}(2,2)\setlength{\unitlength}{0.14in}
\put(0.3,0.3){\circle*{1}}
\end{picture}
$ \ \  denotes a non-isotropic odd simple root (cf.~\cite{Kac77}). We use the same notations for the associated data as in the case for $\mf{g}_\infty$ (${\mf g}={\mf b}, {\mf c}, {\mf d}$).

For $\mf{g}={\mf b}, {\mf b}^\bullet, {\mf c}, {\mf d}$ and $k\in\N$, let ${\mf{g}}_k$ be the subalgebra of ${\mf{g}}_\infty$ whose Dynkin diagram  corresponds to the simple roots $\alpha_0,\ldots,\alpha_{k-1}$.
We assume that $k\geq 2$ for ${\mf g}={\mf b}, {\mf b}^\bullet, {\mf c}$ and $k\geq 4$ for ${\mf g}={\mf d}$.
Let $\mf{l}_k=\mf{l}_\infty\cap\mf{g}_k$ be the corresponding subalgebra of type $A_{k-1}$.

Let $V$ and $W$ be modules over a Lie algebra $\mf g$, where $W$ is irreducible. We define $[V:W]=\dim {\rm Hom}_{\mf g}(W,V)$, the multiplicity of $W$ in $V$.

\subsection{Dual pairs}
Throughout the paper, ${\rm G}_n$ denotes one of the following complex reductive algebraic groups: ${\rm Sp}_n$, ${\rm O}_n$, ${\rm Spin}_n$ and ${\rm Pin}_n$ for $n \geq 2$, where $n$ is even for ${\rm Sp}_n$ and ${\rm Pin}_n$, and $n$ is odd for ${\rm Spin}_n$.  Following \cite[Section 2.2]{CK}, let $V^\la_{{\rm G}_n}$ denote the finite-dimensional irreducible representation of ${\rm G}_n$ corresponding to $\la\in {\mc P}({\rm G}_n)$ (see also \cite{GW,H}), where
\begin{equation*}\label{eq:P(G)}
\begin{split}
{\mc P}({\rm Sp}_n)&=\{\,\la=(\la_1,
\cdots,\la_{\frac{n}{2}})\,|\,\lambda_i\in\Z_+, \ \la_1\geq\ldots \geq
\la_{\frac{n}{2}}\,\}, \\
{\mc P}({\rm Pin}_n)&=\{\,\la=(\la_1,
\cdots,\la_{\frac{n}{2}})\,|\,\lambda_i\in\Z_+, \ \la_1\geq\ldots \geq
\la_{\frac{n}{2}}\,\},\\
{\mc P}({\rm Spin}_n)&=\{\,\la=(\la_1,
\cdots,\la_{\frac{n-1}{2}})\,|\,\lambda_i\in\Z_+, \ \la_1\geq\ldots \geq
\la_{\frac{n-1}{2}}\,\},\\
{\mc P}({\rm O}_n)&=\{\,\la=(\la_1,
\cdots,\la_n)\,|\,\lambda_i\in\Z_+, \ \la_1\geq\ldots \geq \la_{n},\
\la'_1+\la'_2\leq n\,\}.
\end{split}
\end{equation*}
Here $\la'=(\la'_1,\la'_2,\ldots)$ is the conjugate partition of $\la=(\la_1,\la_2,\ldots)$.
From now on, we mean by $({\mf g}_\infty,{\rm G}_n)$ one of the pairs
\begin{equation}\label{eq:dual pairs}
({\mf c}_\infty, {\rm Sp}_n),\ \ ({\mf b}_\infty, {\rm Pin}_n),\ \ ({\mf b}_\infty, {\rm Spin}_n),\ \  ({\mf d}_\infty, {\rm O}_n),
\end{equation}
unless otherwise specified.
For $\la\in {\mc P}({\rm G}_n)$, we define a dominant integral weight $\La^{\mf g}(\la)$ for ${\mf g}_\infty$ by
\begin{equation}\label{eq:hw}
\Lambda^{\mf g}(\la)=
\frac{n}{\epsilon}\La^{\mf g}_0+\la'_1\epsilon_1+\la'_2\epsilon_2+\cdots  ,
\end{equation}
where $\epsilon=2$ if ${\mf g}={\mf c}$, and $\epsilon=1$ otherwise. Note that $\bigsqcup_{n}{\mc P}({\rm G}_n)$ parameterizes the set of all dominant integral weights for ${\mf g}_\infty$. We denote by $L({\mf g}_\infty, \La^{\mf g}(\la))$ the irreducible highest weight ${\mf g}_\infty$-module with highest weight $\La^{\mf g}(\la)$.

For a positive integer $\ell\geq 1$, let
\begin{align*}
\psi^{+,i}(z)&=\sum_{r\in\hf+\Z}\psi^{+,i}_rz^{-r-\hf},\quad
\psi^{-,i}(z)=\sum_{s\in\hf+\Z}\psi^{-,i}_sz^{-s-\hf}\quad (1\leq i\leq \ell)
\end{align*}
be $\ell$ pairs of free fermions
with non-trivial commutation relations
$[\psi^{+,i}_r,\psi^{-,j}_s]=\delta_{ij}\delta_{r+s,0}$. Let
$\mc{F}^\ell$ denote the corresponding Fock space generated by the
vacuum vector $|0\rangle$, which is annihilated by
$\psi^{+,i}_r,\psi^{-,i}_s$ for $r,s>0$.
We introduce a neutral fermionic field
$\phi(z)=\sum_{r\in\hf+\Z}\phi_rz^{-r-\hf}$ with non-trivial commutation
relations $[\phi_r,\phi_s]=\delta_{r+s,0}$. Denote by $\mc{F}^{\hf}$ the
Fock space of $\phi(z)$ generated by a vacuum vector that is
annihilated by $\phi_r$ for $r> 0$. We denote by $\mc{F}^{\ell+\hf}$
the tensor product of $\mc{F}^{\ell}$ and $\mc{F}^{\hf}$. Then we have the following $({\mf g}_\infty,{\rm G}_n)$-duality on $\mc{F}^{\frac{n}{2}}$.

\begin{prop} \label{duality} \cite{Wa}  There exists an action of
$\mf{g}_\infty\times {\rm G}_n$ on $\mc{F}^{\frac{n}{2}}$.
Furthermore, under this joint action, we have
\begin{equation*} \label{eq:dual}
\mc{F}^{\frac{n}{2}}\cong\bigoplus_{\la\in {\mc P}({\rm G}_n)}
L(\mf{g}_\infty,\La^{\mf g}(\la))\otimes V_{{\rm G}_n}^\la.
\end{equation*}
\end{prop}

\subsection{Crystals}
Let us give a brief review on crystals (see \cite{HK,Kas94} for more details).
Let $\mathfrak{g}$ be a
Kac-Moody algebra associated to a symmetrizable generalized Cartan matrix
$A=(a_{ij})_{i,j\in I}$ indexed by a set $I$. Let $P^\vee$ be  the dual weight lattice,
$P={\rm Hom}_\Z(P^\vee,\Z)$ the weight lattice,
$\Pi^\vee=\{\,\alpha^\vee_i\,|\,i\in I\,\}$ the set of simple coroots, and
$\Pi=\{\,\alpha_i\,|\,i\in I\,\}$ the set of simple roots of $\g$ such that $\langle\alpha_j,\alpha^\vee_i\rangle=a_{ij}$ for $i,j\in I$. Let $U_q(\g)$ be the quantized enveloping algebra of $\g$.

A {\it $\g$-crystal} (or {\it crystal} for short) is a set
$B$ together with the maps ${\rm wt} : B \rightarrow P$,
$\varepsilon_i, \varphi_i: B \rightarrow \mathbb{Z}\cup\{-\infty\}$ and
$\te_i, \tf_i: B \rightarrow B\cup\{{\bf 0}\}$ ($i\in I$) satisfying certain axioms.
For a dominant integral weight $\Lambda$ for $\g$, we denote by ${\bf B}(\g,\Lambda)$ the $\g$-crystal  associated to an irreducible highest weight $U_q(\g)$-module with highest weight
$\Lambda$.

Let $B_1$ and $B_2$ be crystals.
A { tensor product $B_1\otimes B_2$}
is defined to be $B_1\times B_2$  as a set with elements  denoted by
$b_1\otimes b_2$, where  {\allowdisplaybreaks
\begin{equation*}
\begin{split}
{\rm wt}(b_1\otimes b_2)&={\rm wt}(b_1)+{\rm wt}(b_2), \\
\varepsilon_i(b_1\otimes b_2)&= {\rm
max}\{\varepsilon_i(b_1),\varepsilon_i(b_2)-\langle {\rm
wt}(b_1),\alpha^\vee_i\rangle\}, \\
\varphi_i(b_1\otimes b_2)&= {\rm max}\{\varphi_i(b_1)+\langle {\rm
wt}(b_2),\alpha^\vee_i\rangle,\varphi_i(b_2)\},\\
{\te}_i(b_1\otimes b_2)&=
\begin{cases}
{\te}_i b_1 \otimes b_2, & \text{if $\varphi_i(b_1)\geq \varepsilon_i(b_2)$}, \\
b_1\otimes {\te}_i b_2, & \text{if
$\varphi_i(b_1)<\varepsilon_i(b_2)$},
\end{cases}\\
{\tf}_i(b_1\otimes b_2)&=
\begin{cases}
{\tf}_i b_1 \otimes b_2, & \text{if  $\varphi_i(b_1)>\varepsilon_i(b_2)$}, \\
b_1\otimes {\tf}_i b_2, & \text{if $\varphi_i(b_1)\leq
\varepsilon_i(b_2)$},
\end{cases}
\end{split}
\end{equation*}
\noindent for $i\in I$. Here we assume that ${\bf 0}\otimes
b_2=b_1\otimes {\bf 0}={\bf 0}$.} Then $B_1\otimes B_2$ is a
crystal.

For  $b_i\in B_i $ ($i=1,2$), we write $b_1 \equiv b_2$ if there is an isomorphism of crystals $C(b_1) \rightarrow C(b_2)$ mapping $b_1$ to $b_2$, where $C(b_i)$ denotes the connected component of $b_i$ in $B_i$.\vskip 2mm

\subsection{Crystals of Young tableaux and Littlewood-Richardson coefficients}\label{subsec:tableaux}

Let $\cP$ be the set of partitions. For $\lambda=(\lambda_i)_{i\geq 1}\in \cP$, we denote by $\ell(\lambda)$ the length of $\lambda$. We identify $\lambda$ with a Young diagram as usual. For a skew Young diagram $\lambda/\mu$, let $SST(\lambda/\mu)$ be the set of semistandard tableaux of shape $\lambda/\mu$ with entries in $\N$.
For $T\in SST(\lambda/\mu)$, let 
$w(T)$ be the column word of $T$, that is, the word given by reading the entries of $T$ column by column from right to left and from top to bottom in each column.  For $T\in SST(\lambda)$ and
$a\in \N$, we denote by $a \rightarrow T$ the tableau obtained by the column insertion of $a$ into $T$ (cf. \cite{Ful}). For a semistandard tableau $S$, we define $(S\rightarrow T)=(w(S)\rightarrow T)$.

Considering $\N$ as a ${\mf l}_{\infty}$-crystal (of type $A_{+\infty}$) associated to the natural representation of ${\mf l}_\infty$, we may regard $SST(\lambda/\mu)$ as an ${\mf l}_{\infty}$-crystal \cite{KashNaka}. In particular, $SST(\lambda)$  is isomorphic to $\B({\mf l}_\infty, \la)$  for $\la\in \cP$, where we identify $\la=(\la_i)_{i\geq 1}$ with $\sum_{i\geq 1}\la_i\epsilon_i$. We denote by $H_\la$ the highest weight element in $SST(\la)$ with weight $\sum_{i\geq 1}\la_i\epsilon_i$, where each $i$th row is filled with $i$.
For $T\in SST(\la)$ and a semistandard tableau $S$, we have $(S\rightarrow T)\equiv T\otimes S$.

Let $\texttt{LR}^\la_{\mu \nu}$ be the set of Littlewood-Richardson tableaux corresponding to $\la, \mu, \nu \in \cP$, which is the set of tableaux $T \in SST(\la/\mu)$ with weight $\nu$ such that $w(T)$ is a lattice word (cf.~\cite{Ful}). Let $c^\la_{\mu \nu}=|\texttt{LR}^\la_{\mu \nu}|$ be the Littlewood-Richardson coefficient corresponding to $\la, \mu, \nu \in \cP$. Recall that $c^\la_{\mu\nu}$ is the number of connected components in $SST(\mu)\otimes SST(\nu)$, which is isomorphic to $SST(\la)$ as an ${\mf l}_\infty$-crystal. By tensor product rule of crystal, we may regard $\texttt{LR}^\la_{\mu \nu}$ as the set of $T\in SST(\nu)$ such that   ${\rm wt}(H_\mu) + {\rm wt}(T)={\rm wt}(H_\la)$, and $\varepsilon_i(T)\leq \varphi_i(H_\mu)=\mu_i-\mu_{i+1}$ for all $i\geq 1$  \cite{Na}.

\section{Spinor model for crystals of types $B, C, D$}

In this section, we briefly review the combinatorial model for ${\bf B}({\mf g}_\infty,\La^{\mf g}(\la))$ ($\la\in \mc{P}({\rm G}_n)$), which was introduced by the author in \cite{K13,K14} (see also \cite[Section 2]{Le} for more details on the existence of the crystal ${\bf B}({\mf g}_\infty,\La^{\mf g}(\la))$ for $U_q(\mf g_\infty)$).

\subsection{Crystals of fundamental representations of ${\mf g}_\infty$}
In this subsection, we describe a combinatorial model for $\B({\mf g}_\infty,\La^{\mf g}_a)$   ($a\in \Z_+$). It is given by translating the $q$-deformed Fock space models for classical Lie algebras due to Hayashi \cite{Ha} in terms of semistandard tableaux of skew shapes with at most two columns.

For a single-columned tableau $U$, let ${\rm ht}(U)$ denote the height of $U$, and $U(i)$ (resp. $U[i]$) the $i$th entry of $U$ from the bottom (resp. the top) for $i\geq 1$.

For $a,b,c\in\Z_+$, let $\lambda(a,b,c)=(2^{b+c},1^a)/(1^b)$, a skew Young diagram with heights $a+c$ and $b+c$ from the left, where  $a$ and $b$ denote the heights of lower and upper single columns in the diagram, respectively.  For example,

$$
\lambda(2,1,3)\ =\
\resizebox{.055\hsize}{!}{$
{\def\lr#1{\multicolumn{1}{|@{\hspace{.75ex}}c@{\hspace{.75ex}}|}{\raisebox{-.04ex}{$#1$}}}\raisebox{-.6ex}
{$\begin{array}{cccc}
\cline{2-2}
\cdot &\lr{}\\
\cline{1-1}\cline{2-2}
\lr{} &\lr{}\\
\cline{1-1}\cline{2-2}
\lr{} &\lr{\ \ }\\
\cline{1-1}\cline{2-2}
\lr{}&\lr{}\\
\cline{1-1}\cline{2-2}
\lr{\ \ }& \cdot \\
\cline{1-1}
\lr{} & \cdot \\
\cline{1-1}
\end{array}$}}$}$$

Let  $T$ be a tableau of shape $\lambda(a,b,c)$, each of whose column is a semistandard tableau. We denote by $T^{\tt L}$ and $T^{\tt R}$ the left and right columns of $T$, respectively.  We also denote by $T^{\tt tail}$ the subtableau of $T$ corresponding to the tail of $\la(a,b,c)$, a lower single column of height $a$, and denote by $T^{\tt body}$ the subtableau of $T$ above $T^{\tt tail}$.

$$T=\resizebox{.05\hsize}{!}
{\def\lr#1{\multicolumn{1}{|@{\hspace{.75ex}}c@{\hspace{.75ex}}|}{\raisebox{-.04ex}{$#1$}}}\raisebox{-.6ex}
{$\begin{array}{cc}
\cline{2-2}
&\lr{3}\\
\cline{1-1}\cline{2-2}
\lr{2} &\lr{4}\\
\cline{1-1}\cline{2-2}
\lr{4}&\lr{6}\\
\cline{1-1}\cline{2-2}
\lr{5} & \lr{8}\\
\cline{1-1}\cline{2-2}
\lr{6} \\
\cline{1-1}
\lr{7} \\
\cline{1-1}
\end{array}$}}\ \ \ \ \ \ \ \ \
T^{\tt body}=\resizebox{.05\hsize}{!}
{\def\lr#1{\multicolumn{1}{|@{\hspace{.75ex}}c@{\hspace{.75ex}}|}{\raisebox{-.04ex}{$#1$}}}\raisebox{-.6ex}
{$\begin{array}{cc}
\cline{2-2}
&\lr{3}\\
\cline{1-1}\cline{2-2}
\lr{2} &\lr{4}\\
\cline{1-1}\cline{2-2}
\lr{4}&\lr{6}\\
\cline{1-1}\cline{2-2}
\lr{5} & \lr{8}\\
\cline{1-1}\cline{2-2}
\end{array}$}}\ \ \ \ \
T^{\tt tail}=\ \resizebox{.025\hsize}{!}
{\def\lr#1{\multicolumn{1}{|@{\hspace{.75ex}}c@{\hspace{.75ex}}|}{\raisebox{-.04ex}{$#1$}}}\raisebox{-.6ex}
{$\begin{array}{c}
\cline{1-1}
\lr{6} \\
\cline{1-1}
\lr{7} \\
\cline{1-1}
\end{array}$}}\ \ \ \ \
$$\vskip 3mm

Suppose that $T$ is semistandard, that is, $T\in SST(\la(a,b,c))$. One may slide down $T^{\tt R}$ by $k$ positions for $0\leq k\leq \min\{a,b\}$ to have a (not necessarily semistandard) tableau $T'$ of shape $\lambda(a-k,b-k,c+k)$. We define the {\it residue of $T$} to be the maximal $k$ such that
$T'$ is semistandard, and denote it by ${\mf r}_T$.
For example, ${\mf r}_{S}=1$ and ${\mf r}_T=2$, when

$$S=\resizebox{.05\hsize}{!}
{\def\lr#1{\multicolumn{1}{|@{\hspace{.75ex}}c@{\hspace{.75ex}}|}{\raisebox{-.04ex}{$#1$}}}\raisebox{-.6ex}
{$\begin{array}{cc}
\cline{2-2}
&\lr{3}\\
 \cline{2-2}
 \cline{2-2}
 &\lr{4}\\
\cline{1-1}\cline{2-2}
\lr{4}&\lr{6}\\
\cline{1-1}\cline{2-2}
\lr{5} & \lr{8}\\
\cline{1-1}\cline{2-2}
\lr{6} \\
\cline{1-1}
\lr{7} \\
\cline{1-1}
\end{array}$}}\ \ \ \ \ \ \ \ \
T=\resizebox{.05\hsize}{!}
{\def\lr#1{\multicolumn{1}{|@{\hspace{.75ex}}c@{\hspace{.75ex}}|}{\raisebox{-.04ex}{$#1$}}}\raisebox{-.6ex}
{$\begin{array}{cc}
\cline{2-2}
&\lr{4}\\
 \cline{2-2}
 \cline{2-2}
 &\lr{6}\\
\cline{1-1}\cline{2-2}
\lr{4}&\lr{7}\\
\cline{1-1}\cline{2-2}
\lr{5} & \lr{8}\\
\cline{1-1}\cline{2-2}
\lr{6} \\
\cline{1-1}
\lr{7} \\
\cline{1-1}
\end{array}$}}
$$\vskip 3mm

For  $a \in \Z_+$, let
\begin{equation*}
\begin{split}
&{\bf T}^{{\mf g}}(a)=\{\,T\,|\,T\in SST(\la(a,b,c)), \  (b,c)\in {\mc H}^{\mf g}, \  {\mf r}_T\leq r^{\mf g}   \,\},\\
\end{split}
\end{equation*}
where 
\begin{equation*}
{\mc H}^{\mf g}=
\begin{cases}
\{0\}\times  \Z_+ , & \text{if ${\mf g}={\mf c}$},\\
\Z_+\times  \Z_+ , & \text{if ${\mf g}={\mf b}$},\\
2\Z_+\times  2\Z_+ , & \text{if ${\mf g}={\mf d}$},\\
\end{cases}\ \ \ \ \ \
r^{\mf g}=
\begin{cases}
0 , & \text{if ${\mf g}={\mf c}, {\mf b}$},\\
1 , & \text{if ${\mf g}={\mf d}$}.
\end{cases}
\end{equation*}
We also put
\begin{equation*}
\begin{split}
&{\bf T}^{\rm sp}= \bigsqcup_{a\,\in \Z_+} SST((1^a)), \\ & {\bf T}^{\rm sp\,+} =\{\,T\in {\bf T}^{\rm sp}\,|\,{\mf r}_{T}=0\,\}, \ {\bf T}^{\rm sp\,-}=\{\,T\in {\bf T}^{\rm sp}\,|\,{\mf r}_{T}=1\,\},\\
&\ov{\bf T}^{\mf d}(0)=\bigsqcup_{(b,c)\in {\mc H}^{\mf d}}SST(\la(0,b,c+1)),\\
\end{split}
\end{equation*}
Here we define the residue ${\mf r}_T$ of $T\in {\bf T}^{\rm sp}$ to be the residue of ${\rm ht}(T)$ modulo 2. (Note that if $T\in SST(\la(0,b,0))\subset T^{\g}(0)$, then the shape of $T$ is a single column but ${\mf r}_T=0$.)

Now, let ${\bf B}$ be one of ${\bf T}^{{\mf g}}(a)$ ($a\in\Z_+$), ${\bf T}^{\rm sp}$, and $\ov{\bf T}^{\mf d}(0)$, and let $T\in {\bf B}$ be given. 
First, we define $\te_i T$ and $\tf_i T$ for $i\in I\setminus\{0\}$ regarding ${\bf B}$ as a subset of an ${\mf l}_\infty$-crystal $\bigsqcup_{\la\in\cP}SST(\la)$. Next, we define $\te_0T$ and $\tf_0T$ as follows:

\begin{itemize}
\item[(1)]  Suppose that ${\mf g}={\mf c}$ and ${\bf B}={\bf T}^{\mf c}(a)$.  We define $\te_0 T$ to be the tableau obtained from $T$ by removing a domino
$\resizebox{.05\hsize}{!}{\def\lr#1{\multicolumn{1}{|@{\hspace{.6ex}}c@{\hspace{.6ex}}|}{\raisebox{-.25ex}{$#1$}}}\raisebox{-.65ex}
{$\begin{array}[b]{cc}
\cline{1-1}\cline{2-2}
\lr{1}&\lr{1}\\
\cline{1-1}\cline{2-2}
\end{array}$}}$\ if $T$ has $\resizebox{.05\hsize}{!}{\def\lr#1{\multicolumn{1}{|@{\hspace{.6ex}}c@{\hspace{.6ex}}|}{\raisebox{-.25ex}{$#1$}}}\raisebox{-.65ex}
{$\begin{array}[b]{cc}
\cline{1-1}\cline{2-2}
\lr{1}&\lr{1}\\
\cline{1-1}\cline{2-2}
\end{array}$}}$ on its top, and ${\bf 0}$, otherwise.
We define $\tf_0T$ in a similar way by adding $\resizebox{.05\hsize}{!}{\def\lr#1{\multicolumn{1}{|@{\hspace{.6ex}}c@{\hspace{.6ex}}|}{\raisebox{-.25ex}{$#1$}}}\raisebox{-.65ex}
{$\begin{array}[b]{cc}
\cline{1-1}\cline{2-2}
\lr{1}&\lr{1}\\
\cline{1-1}\cline{2-2}
\end{array}$}}$ on top of $T$.

\item[(2)] Suppose that ${\mf g}={\mf b}$. When ${\bf B}={\bf T}^{\rm sp}$, we define
 $\te_0 T$ to be the tableau obtained from $T$ by removing  \boxed{$1$} \ if $T$ has \boxed{$1$} on its top, and ${\bf 0}$, otherwise. 
We define $\tf_0T$ in a similar way by adding \boxed{$1$} on top of $T$. 
When ${\bf B}={\bf T}^{\mf b}(a)$, we define $\tilde{x}_0 T$ ($x=e, f$) to be $\tilde{x}_0(T^{\tt R}\otimes T^{\tt L})$ regarding ${\bf T}^{\rm sp}$ as a regular $\mf{sl}_2$-crystal with respect to $\te_0$ and $\tf_0$, and applying the tensor product rule of crystals to ${\bf B}\subset ({\bf T}^{\rm sp})^{\otimes 2}$.

\item[(3)] Suppose that ${\mf g}={\mf d}$. When ${\bf B}={\bf T}^{\rm sp}$, we define
$\te_0 T$ to be the tableau obtained from $T$ by removing a domino
$
\resizebox{.02\hsize}{!}{$
{\def\lr#1{\multicolumn{1}{|@{\hspace{.75ex}}c@{\hspace{.75ex}}|}{\raisebox{-.04ex}{$#1$}}}\raisebox{0.5ex}
{$\begin{array}{cc}
\cline{1-1}
\lr{1} \\
\cline{1-1}
\lr{2} \\
\cline{1-1}
\end{array}$}}$}
$\, if $T$ has
$
\resizebox{.02\hsize}{!}{$
{\def\lr#1{\multicolumn{1}{|@{\hspace{.75ex}}c@{\hspace{.75ex}}|}{\raisebox{-.04ex}{$#1$}}}\raisebox{0.5ex}
{$\begin{array}{cc}
\cline{1-1}
\lr{1} \\
\cline{1-1}
\lr{2} \\
\cline{1-1}
\end{array}$}}$}
$\, on its top, and ${\bf 0}$ otherwise.
We define $\tf_0T$ in a similar way by adding 
$
\resizebox{.02\hsize}{!}{$
{\def\lr#1{\multicolumn{1}{|@{\hspace{.75ex}}c@{\hspace{.75ex}}|}{\raisebox{-.04ex}{$#1$}}}\raisebox{0.5ex}
{$\begin{array}{cc}
\cline{1-1}
\lr{1} \\
\cline{1-1}
\lr{2} \\
\cline{1-1}
\end{array}$}}$}
$
on top of $T$.
When ${\bf B}={\bf T}^{\mf d}(a)$ or $\ov{{\bf T}}^{\mf d}(0)$, we define $\tilde{x}_0 T$ ($x=e, f$) to be $\tilde{x}_0(T^{\tt R}\otimes T^{\tt L})$ as in (2).
\end{itemize}
We put
\begin{equation*}
{\rm wt}(T)=
\begin{cases}
\frac{2}{\epsilon}\La^{\mf g}_{0}+\sum_{i\geq 1}m_i\epsilon_i, & \text{if  ${\bf B}={\bf T}^{\mf g}(a)$ or $\ov{\bf T}^{\mf d}(0)$},\\
\La^{\mf g}_{0}+\sum_{i\geq 1}m_i\epsilon_i, & \text{if  $\bf{B}={\bf T}^{\rm sp}$},\\
\end{cases}
\end{equation*}
where $m_i$ is the number of occurrences of $i$ in $T$ (cf.~\eqref{eq:hw}), and 
\begin{equation*}
\varepsilon_i(T)=\max\{\,k\,|\,\te_i^kT\neq {\bf 0}\,\},\quad 
\varphi_i(T)=\max\{\,k\,|\,\tf_i^kT\neq {\bf 0}\,\}\quad (i\in I).
\end{equation*}
Then ${\bf B}$ is a $\g_\infty$-crystal with respect to ${\rm wt}$, $\varepsilon_i$, $\varphi_i$, $\te_i$, $\tf_i$ ($i\in I$).   
By \cite[Theorem 7.1]{K13} for ${\mf g}={\mf b}, {\mf c}$ and \cite[Proposition 4.2]{K14} for ${\mf g}={\mf d}$, we have

\begin{prop} \mbox{} 
\begin{itemize}
\item[(1)] For ${\mf g}={\mf c}$, $${\bf T}^{\mf c}(a) \cong \B({\mf c}_\infty, \La^{\mf c}_a)   \ \ \ (a\in \Z_+).$$

\item[(2)] For ${\mf g}={\mf b}$,
\begin{equation*}
\begin{split}
\begin{cases}
{\bf T}^{\mf b}(a)   \cong \B({\mf b}_\infty, \La^{\mf b}_a) \ \ \ (a\geq 1), \\
{\bf T}^{\mf b}(0)   \cong \B({\mf b}_\infty, 2\La^{\mf b}_0),
\end{cases}
\ \
{\bf T}^{\rm sp}   \cong \B({\mf b}_\infty, \La^{\mf b}_0).
\end{split}
\end{equation*}

\item[(3)] For ${\mf g}={\mf d}$,
\begin{equation*}
\begin{split}
\begin{cases}
{\bf T}^{\mf d}(a)   \cong \B({\mf d}_\infty, \La^{\mf d}_a) \ \ \ (a\geq 2), \\
{\bf T}^{\mf d}(0)   \cong \B({\mf d}_\infty, 2\La^{\mf d}_0), \\
\ov{\bf T}^{\mf d}(0)   \cong \B({\mf d}_\infty, 2\La^{\mf d}_1),
\end{cases}
\ \
\begin{cases}
{\bf T}^{\rm sp +}   \cong \B({\mf d}_\infty, \La^{\mf d}_0), \\
{\bf T}^{\rm sp -}   \cong \B({\mf d}_\infty, \La^{\mf d}_1).
\end{cases}
\end{split}
\end{equation*}

\end{itemize}
\end{prop}
\vskip 2mm

Note that the highest weight element in ${\bf B}$ is given by $H_{(1^a)}$, $\emptyset$,\, $H_{(1)}$, and $H_{(2)}$, when ${\bf B}={\bf T}^{\mf g}(a)$, ${\bf T}^{\rm sp,+}$, ${\bf T}^{\rm sp,-}$, and $\ov{\bf T}^{\mf d}(0)$, respectively.

\subsection{Admissibility}
To describe a ${\mf g}_\infty$-crystal  $\B({\mf g}_\infty,\La)$ associated to arbitrary integral dominant weight $\La=\La^{\mf g}_{a_1}+\cdots+\La^{\mf g}_{a_r}$, we embed $\B({\mf g}_\infty,\La)$ into the tensor product  $\B({\mf g}_\infty,\La^{\mf g}_{a_1})\otimes \cdots \otimes \B({\mf g}_\infty,\La^{\mf g}_{a_r})$, and then describe a connected component of the highest weight element with weight $\La$. For this, we need to characterize an explicit condition for $T_i\otimes T_j \in \B({\mf g}_\infty,\La^{\mf g}_{a_i})\otimes \B({\mf g}_\infty,\La^{\mf g}_{a_j})$ to be in the same connected component including the highest weight element of weight $\La^{\mf g}_{a_i}+\La^{\mf g}_{a_j}$. This condition is called {\it admissibility}, which is an analogue of semistandardness between two adjacent columns in type $A$. \vskip 2mm

For $a\in\Z_+$, let $T\in {\bf T}^{\mf g}(a)$ be given.
Associated to $T$, we introduce the following two pairs of single-columned tableaux $({}^{\tt L}T, {}^{\tt R}T)$ and $(T^{\tt L^*}, T^{\tt R^*})$, which will play a crucial role in describing admissibility.

First, we define $({}^{\tt L}T, {}^{\tt R}T)$ to be the pair determined by the following algorithm:
\begin{itemize}
\item[(1)] Let $y_i=T^{\tt R}(i)$ for $1\leq i\leq {\rm ht}(T^{\tt R})$. First, slide down the box $\boxed{y_1}$ in $T^{\tt R}$ as far as the entry of $T^{\tt L}$ in the same row is no greater than $y_1$. If  no entry of $T^{\tt L}$ is greater than $y$, we place $\boxed{y_1}$ next to the bottom of $T^{\tt L}$.

\item[(2)] Next, slide down $\boxed{y_2}$ until it is above $\boxed{y_1}$ and the entry of $T^{\tt L}$ in the same row is no greater than $y_2$. Repeat the same process with the other boxes $\boxed{y_3}, \boxed{y_4},\ldots$  until there is no moving down.

\item[(3)] Slide each box $\boxed{x}$ in $T^{\tt L}$ to the right if its right position is empty  (indeed the number of such boxes is $a-{\mf r}_T$).

\item[(4)] Define ${}^{\tt R}T$ to be the tableau determined by the boxes $\boxed{y_i}$\,'s in $T^{\tt R}$ together with boxes $\boxed{x}$'s which have moved from $T^{\tt L}$ by (3), and ${}^{\tt L}T$ to be the tableau with the remaining boxes on the left.
\end{itemize}
For example, we have 
$$
T\ = \
\resizebox{.43\hsize}{!}{$
{\def\lr#1{\multicolumn{1}{|@{\hspace{.75ex}}c@{\hspace{.75ex}}|}{\raisebox{-.04ex}{$#1$}}}\raisebox{-.6ex}
{$\begin{array}{cc}
\cline{2-2}
&\lr{2}\\
 \cline{2-2}
 &\lr{3}\\
\cline{1-1}\cline{2-2}
\lr{2}&\lr{7}\\
\cline{1-1}\cline{2-2}
\lr{4}&\lr{9}\\
\cline{1-1}\cline{2-2}
\lr{6} \\
\cline{1-1}
\lr{8} \\
\cline{1-1}
\lr{9} \\
\cline{1-1}\\
\!\!\! T^{\tt L}\!\!\! &  T^{\tt R}\!\!\!\!\!
\end{array}$}}
\ \ \  \rightarrow  \ \ \
{\def\lr#1{\multicolumn{1}{|@{\hspace{.75ex}}c@{\hspace{.75ex}}|}{\raisebox{-.04ex}{$#1$}}}\raisebox{-.6ex}
{$\begin{array}{cc}
& \\
\cline{2-2}
&\lr{2}\\
\cline{1-1}\cline{2-2}
\lr{2}& \lr{3}\\
\cline{1-1}\cline{2-2}
\lr{4}\\
\cline{1-1}\cline{2-2}
\lr{6}&\lr{7} \\
\cline{1-1}\cline{2-2}
\lr{8} \\
\cline{1-1}\cline{2-2}
\lr{9}&\lr{9} \\
\cline{1-1}\cline{2-2}\\ \\
\end{array}$}}
\ \ \  \rightarrow  \ \ \
{\def\lr#1{\multicolumn{1}{|@{\hspace{.75ex}}c@{\hspace{.75ex}}|}{\raisebox{-.04ex}{$#1$}}}\raisebox{-.6ex}
{$\begin{array}{cc}
& \\
\cline{2-2}
&\lr{2}\\
\cline{1-1}\cline{2-2}
\lr{2}& \lr{3}\\
\cline{1-1}\cline{2-2}
& \lr{4}\\
\cline{1-1}\cline{2-2}
\lr{6}&\lr{7} \\
\cline{1-1}\cline{2-2}
&\lr{8} \\
\cline{1-1}\cline{2-2}
\lr{9}&\lr{9} \\
\cline{1-1}\cline{2-2}\\ \\
\end{array}$}}
\ \ \  \rightarrow  \ \ \
{\def\lr#1{\multicolumn{1}{|@{\hspace{.75ex}}c@{\hspace{.75ex}}|}{\raisebox{-.04ex}{$#1$}}}\raisebox{-.6ex}
{$\begin{array}{cc}
& \\
\cline{1-1}\cline{2-2}
\lr{2}&\lr{2}\\
\cline{1-1}\cline{2-2}
\lr{6}& \lr{3}\\
\cline{1-1}\cline{2-2}
\lr{9}& \lr{4}\\
\cline{1-1}\cline{2-2}
&\lr{7} \\
\cline{2-2}
&\lr{8} \\
\cline{2-2}
&\lr{9} \\
\cline{2-2}\\
\!\!\! \!{}^{\tt L}T \! &  \!\!{}^{\tt R}T\!\!\!\! 
\end{array}$}}$}
$$
\noindent where the pairs $(T^{\tt L},{}^{\tt R}T)$ and $({}^{\tt L}T, T^{\tt R})$ are arranged to share the same bottom lines, respectively.

Next, when ${\mf r}_T=1$, we define $(T^{\tt L^*}, T^{\tt R^*})$ to be the pair determined by the following algorithm:
\begin{itemize}
\item[(1)] Let $x_i=T^{\tt L}[i]$ for $1\leq i\leq {\rm ht}(T^{\tt L})$.
First, slide upward $\boxed{x_1}$ until the entry of $T^{\tt R}$ in the same row is no smaller than   $x_1$. If  no entry of $T^{\tt L}$ is smaller than $x_1$, we place $\boxed{x_1}$ next to the top of $T^{\tt R}$.

\item[(2)] Next, slide upward $\boxed{x_2}$ until it is below $\boxed{x_1}$ and the entry of $T^{\tt R}$ in the same row is no smaller than $x_2$. Repeat the same process with the other boxes $\boxed{x_3}, \boxed{x_4},\ldots$  until there is no moving up.

\item[(3)] Choose the lowest box $\boxed{y}$ in $T^{\tt R}$ whose left position is empty, and then slide it to the left (there exists at least one such $\boxed{y}$ since ${\mf r}_T=1$).

\item[(4)] Define ${T}^{\tt L^\ast}$ to be the tableau determined by the boxes $\boxed{x_i}$\,'s in $T^{\tt L}$ together with $\boxed{y}$, and define ${T}^{\tt R^\ast}$ to be the tableau given by the remaining boxes on the right.
\end{itemize}

For example, we have \vskip -2mm
$$
T\ = \
\resizebox{.43\hsize}{!}{$
{\def\lr#1{\multicolumn{1}{|@{\hspace{.75ex}}c@{\hspace{.75ex}}|}{\raisebox{-.04ex}{$#1$}}}\raisebox{-.6ex}
{$\begin{array}{cc}
\cline{2-2}
&\lr{2}\\
 \cline{2-2}
 &\lr{3}\\
\cline{1-1}\cline{2-2}
\lr{2}&\lr{7}\\
\cline{1-1}\cline{2-2}
\lr{4}&\lr{9}\\
\cline{1-1}\cline{2-2}
\lr{6} \\
\cline{1-1}
\lr{8} \\
\cline{1-1}
\lr{9} \\
\cline{1-1}\\
\!\!\! T^{\tt L}\!\!\! &  T^{\tt R}\!\!\!\!\!
\end{array}$}}
\ \ \  \rightarrow  \ \ \
{\def\lr#1{\multicolumn{1}{|@{\hspace{.75ex}}c@{\hspace{.75ex}}|}{\raisebox{-.04ex}{$#1$}}}\raisebox{-.6ex}
{$\begin{array}{cc}
\cline{1-1}\cline{2-2}
\lr{2}&\lr{2}\\
\cline{1-1}\cline{2-2}
& \lr{3}\\
\cline{1-1}\cline{2-2}
\lr{4} &\lr{7}\\
\cline{1-1}\cline{2-2}
\lr{6}&\lr{9} \\
\cline{1-1}\cline{2-2}
\lr{8} \\
\cline{1-1}
\lr{9}  \\
\cline{1-1} \\ \\ \\
\end{array}$}}
\ \ \  \rightarrow  \ \ \
{\def\lr#1{\multicolumn{1}{|@{\hspace{.75ex}}c@{\hspace{.75ex}}|}{\raisebox{-.04ex}{$#1$}}}\raisebox{-.6ex}
{$\begin{array}{cc}
\cline{1-1}\cline{2-2}
\lr{2}&\lr{2}\\
\cline{1-1}\cline{2-2}
\lr{3} &\\
\cline{1-1}\cline{2-2}
\lr{4} &\lr{7}\\
\cline{1-1}\cline{2-2}
\lr{6}&\lr{9} \\
\cline{1-1}\cline{2-2}
\lr{8} \\
\cline{1-1}
\lr{9}  \\
\cline{1-1} \\ \\ \\
\end{array}$}}
\ \ \  \rightarrow  \ \ \
{\def\lr#1{\multicolumn{1}{|@{\hspace{.75ex}}c@{\hspace{.75ex}}|}{\raisebox{-.04ex}{$#1$}}}\raisebox{-.6ex}
{$\begin{array}{cc}
& \\
\cline{1-1}\cline{2-2}
\lr{2}&\lr{2}\\
\cline{1-1}\cline{2-2}
\lr{3} &\lr{7}\\
\cline{1-1}\cline{2-2}
\lr{4} &\lr{9}\\
\cline{1-1}\cline{2-2}
\lr{6}& \\
\cline{1-1}
\lr{8} \\
\cline{1-1}
\lr{9}  \\
\cline{1-1}\\
\!\!\! T^{\tt L^*} \!\!\!\!\! &  T^{\tt R^*} \!\!\!\!\!\!\!\!\!
\end{array}$}}$}
$$ \vskip 2mm
\noindent Note that the pairs $(T^{\tt L},T^{\tt L^*})$ and $({T}^{\tt R}, T^{\tt R^*})$ are arranged to share the same bottom lines, respectively.

 We refer the reader to \cite[Section 6]{K13} and \cite[Section 3]{K14} for more details on the well-definedness of $({}^{\tt L}T, {}^{\tt R}T)$ and $(T^{\tt L^*}, T^{\tt R^*})$.\vskip 2mm

\begin{df}\label{admissible}{\rm  Let $a, a'\in \Z_+$ be given with $a\geq a'$. We say that a pair $(T,S)$ is admissible, and write  $T\prec S$ if  it is one of the following cases:

\begin{itemize}
\item[(1)] $(T,S)\in {\bf T}^{\mf g}(a)\times {\bf T}^{\mf g}(a')$ or ${\bf T}^{\mf g}(a)\times  {\bf T}^{\rm sp}$ with
\begin{equation*}
\begin{split}
& \ \ \ {\rm ht}({T}^{\tt R})\leq {\rm ht}(S^{\tt L}) -a'+  2{\mf r}_{T}{\mf r}_{S}, \\
& \begin{cases}
{T}^{\tt R}(i)\leq  {}^{\tt L}S(i), & \text{if ${\mf r}_T=0$ or ${\mf r}_S=0$},\\
{T}^{\tt R^\ast}(i)\leq  {}^{\tt L}S(i), & \text{if ${\mf r}_T={\mf r}_S=1$},\\
\end{cases}\\
&\begin{cases}
{}^{\tt R}T(i+a-a')\leq {S}^{\tt L}(i), & \text{if ${\mf r}_T=0$ or ${\mf r}_S=0$},\\
{}^{\tt R}T(i+a-a'+\varepsilon)\leq {S}^{\tt L^\ast}(i), & \text{if ${\mf r}_T={\mf r}_S=1$},\\
\end{cases}
\end{split}
\end{equation*}
for $i\geq 1$ (Here $\varepsilon=1$ if $S\in {\bf T}^{\rm sp -}$  and $0$ otherwise, and we assume that $a'={\mf r}_S$, $S=S^{\tt L}={}^{\tt L}S={S}^{\tt L^\ast}$ when $S\in {\bf T}^{\rm sp}$),\vskip 2mm

\item[(2)]
$(T,S) \in {\bf T}^{\mf{d}}(a)\times \ov{\bf T}^{\mf{d}}(0)$ with $T\prec S^{\tt L}$ in the sense of (1), regarding $S^{\tt L}\in {\bf T}^{\rm sp -}$,
\vskip 2mm

\item[(3)]
$(T,S)\in \ov{\bf T}^{\mf{d}}(0)\times \ov{\bf T}^{\mf{d}}(0)$ or $ \ov{\bf T}^{\mf{d}}(0)\times   {\bf T}^{\rm sp -}$ with $(T^{\tt R},S^{\tt L})\in \ov{\bf T}^{\mf{d}}(0)$.\vskip 2mm
\end{itemize}

}\end{df}

\begin{ex}\label{ex:admissibility}
{\rm
For $T\in {\bf T}^{\mf c}(2)$ and $S\in {\bf T}^{\mf c}(1)$ below, we have
$$\resizebox{.4\hsize}{!}
{\def\lr#1{\multicolumn{1}{|@{\hspace{.75ex}}c@{\hspace{.75ex}}|}{\raisebox{-.04ex}{$#1$}}}\raisebox{-.6ex}
{$\begin{array}{cccccccccccccccccc}
& T\!\!\!\!\!\!\!\!\! & & & & & & & & & & \!\! S\!\!\!\!\!\!\!\!\! & & & & & & \\
\cline{12-12}\cline{13-13}\cline{16-16}\cline{17-17}
&  & & & & & & & & & & \lr{\color{blue} 1} &\lr{3} & & & \lr{\color{red} 2} & \lr{1} \\
\cline{2-2}\cline{3-3}\cline{6-6}\cline{7-7}\cline{12-12}\cline{13-13}\cline{16-16}\cline{17-17}
& \lr{2}&\lr{\color{red} 4} & & & \lr{4} &  \lr{\color{blue} 2} & & & & & \lr{\color{blue}2}&\lr{6} & & & \lr{\color{red}5} &  \lr{3}\\
\cline{2-2}\cline{3-3}\cline{6-6}\cline{7-7}\cline{12-12}\cline{13-13}\cline{16-16}\cline{17-17}
& \lr{3} & \lr{\color{red}6}& & & \lr{5} & \lr{\color{blue}3} & & & & & \lr{\color{blue}5} & \lr{7}& & & \lr{\color{red}6} &  \lr{6}\\
\cdashline{1-1}[0.5pt/1pt]\cline{2-2}\cline{3-3}\cdashline{4-5}[0.5pt/1pt]\cline{6-6}\cline{7-7}\cdashline{8-11}[0.5pt/1pt]
\cline{12-12}\cline{13-13}\cdashline{14-18}[0.5pt/1pt]\cline{16-16}\cline{17-17}
& \lr{4} & & & & & \lr{\color{blue}4} & & & & & \lr{\color{blue}6} & & & & & \lr{7}& \\
\cline{2-2}\cline{7-7}\cline{12-12}\cline{17-17}
& \lr{5} & & & & & \lr{\color{blue}6}& \\
\cline{2-2}\cline{7-7}\\
& \!\!\!\!\!\! T^{\tt L}\!\!\!\!\! & \!\! T^{\tt R}\!\!\!\!\!\! & & & \!\!\!\!\!\! {}^{\tt L}T\!\!\!\!\! & \!\! {}^{\tt R}T\!\!\!\!\!\! & & & & & \!\!\!\!\!\! S^{\tt L}\!\!\!\!\! & \!\! S^{\tt R}\!\!\!\!\!\! & & & \!\!\!\!\!\! {}^{\tt L}S\!\!\!\!\! & \!\! {}^{\tt R}S\!\!\!\!\!\! & \\
\end{array}$}} 
$$ 
\noindent where the dashed line denotes the line separating the body and tail of $T$ and $S$.
Then 
$({}^{\tt R}T,S^{\tt L})$ (in blue) and $(T^{\tt R},{}^{\tt L}S)$  (in red) form semistandard tableaux \vskip 2mm
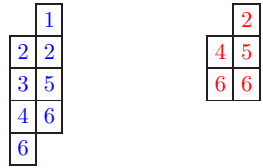
$$
\resizebox{.05\hsize}{!}
{\def\lr#1{\multicolumn{1}{|@{\hspace{.75ex}}c@{\hspace{.75ex}}|}{\raisebox{-.04ex}{$#1$}}}\raisebox{-.6ex}
{$\begin{array}{cc}
 \cline{2-2}
  & \lr{\color{blue}1} \\
\cline{1-1}\cline{2-2}
  \lr{\color{blue}2} & \lr{\color{blue}2}\\
\cline{1-1}\cline{2-2}
 \lr{\color{blue}3}& \lr{\color{blue}5}  \\
\cline{1-1}\cline{2-2}
\lr{\color{blue}4} & \lr{\color{blue}6} \\
\cline{1-1}\cline{2-2}
\lr{\color{blue}6} \\ 
\cline{1-1}
\end{array}$}}
 \ \ \ \ \  \ \ \ \ \  \ \ \ \ \ 
\resizebox{.05\hsize}{!}
{\def\lr#1{\multicolumn{1}{|@{\hspace{.75ex}}c@{\hspace{.75ex}}|}{\raisebox{-.04ex}{$#1$}}}\raisebox{-.6ex}
{$\begin{array}{cc}
 \cline{2-2}
  & \lr{\color{red}2} \\
\cline{1-1}\cline{2-2}
\lr{\color{red}{4}} &  \lr{\color{red}5}\\
\cline{1-1}\cline{2-2}
\lr{\color{red}6} & \lr{\color{red}6}  \\
\cline{1-1}\cline{2-2}\\
\\
\end{array}$}}
$$\vskip 2mm
\noindent
which implies that $T\prec S$.
}
\end{ex}

The following lemma will be used in the next section.
\begin{lem}\label{lem:aux-1}
Let $(T,S)\in {\bf T}^{\mf g}(a)\times {\bf T}^{\mf g}(a')$ $(a\geq a')$ such that $T\prec S$.  Then
\begin{itemize}
\item[(1)] $T^{\tt L}(i+a-a')\leq S^{\tt L}(i)$ for $i\geq 1$, if ${\mf r}_T=0$ or ${\mf r}_S=0$,

\item[(2)] $T^{\tt R}(1)\leq S^{\tt R}(1)$.
\end{itemize}
\end{lem}
\pf (1) By definition of ${}^{\tt R}T$, we have $T^{\tt L}(i+a-a')\leq {}^{\tt R}T(i+a-a')$ for $i\geq 1$. Since ${\mf r}_T=0$ or ${\mf r}_S=0$, we have $T^{\tt L}(i+a-a')\leq S^{\tt L}(i)$ for $i\geq 1$ by Definition \ref{admissible}(1).

(2) Suppose first that ${\mf g}={\mf b}$ or ${\mf c}$. In this case, we always have ${\mf r}_T={\mf r}_S=0$. We have ${}^{\tt L}S(1)\leq S^{\tt R}(1)$ by definition of ${}^{\tt L}S$, and hence $T^{\tt R}(1)\leq {}^{\tt L}S(1)\leq S^{\tt R}(1)$ by Definition \ref{admissible}(1).

Now, we suppose that ${\mf g}={\mf d}$. If  ${\mf r}_T=0$ or ${\mf r}_S=0$, then we also have $T^{\tt R}(1)\leq S^{\tt R}(1)$ by the same argument as in ${\mf g}={\mf b}$ or ${\mf c}$. So we may assume that ${\mf r}_T={\mf r}_S=1$. Let $x=T^{\tt R}(1)$. Suppose that $\boxed{x}$ moves to the right when we construct $T^{\tt L^\ast}$. This implies that $T^{\tt L}(a)<x<T^{\tt L}(a-1)$, and hence ${}^{\tt R}T(a)=x$. On the other hand, we have  $S^{\tt L^\ast}(a')=S^{\tt L}(a')\leq S^{\tt R}(1)$ or $S^{\tt L^\ast}(a')= S^{\tt R}(1)$ since ${\mf r}_S=1$ depending on whether $\boxed{y}$ with $y=S^{\tt R}(1)$ moves to the left or not when we construct $S^{\tt L^\ast}$. By Definition \ref{admissible}(1), we have $x={}^{\tt R}T(a)\leq S^{\tt L^\ast}(a')\leq S^{\tt R}(1)$. Next, suppose that $\boxed{x}$ does not  move to the right when we construct $T^{\tt L^\ast}$. This implies that $T^{\tt R^\ast}(1)=x$. By Definition \ref{admissible}(1), we have $x\leq {}^{\tt L}S(1)\leq S^{\tt R}(1)$. This completes the proof.
\qed

\subsection{Highest weight crystals}
Let $({\mf g}_\infty,{\rm G}_n)$ be one of the dual pairs in \eqref{eq:dual pairs} and $\la\in {\mc P}({\rm G}_n)$.

Suppose that $({\mf g}_\infty,{\rm G}_n)\neq({\mf d}_\infty, {\rm O}_n)$. We put
\begin{equation*}
\widehat{\bf T}^{\mf g}(\la,n)=
\begin{cases}
{\bf T}^{\mf g}(\la_1)\times \cdots \times {\bf T}^{\mf g}(\la_{\frac{n}{2}}), & \text{for $({\mf c}_\infty, {\rm Sp}_n)$ or $({\mf b}_\infty, {\rm Pin}_n)$}, \\
{\bf T}^{\mf g}(\la_1)\times \cdots \times {\bf T}^{\mf g}(\la_{\frac{n-1}{2}})\times {\bf T}^{\rm sp}, & \text{for $({\mf b}_\infty, {\rm Spin}_n)$}. \\
\end{cases}
\end{equation*}

Suppose that $({\mf g}_\infty,{\rm G}_n)=({\mf d}_\infty, {\rm O}_n)$. Let $q_\pm$ and $r_\pm$ be non-negative integers such that
\begin{equation*}
\begin{cases}
n-2\lambda'_1=2q_++r_+, & \text{if $n-2\lambda'_1\geq 0$,}\\
2\lambda'_1-n=2q_-+r_-, & \text{if $n-2\lambda'_1\leq 0$,}\\
\end{cases}
\end{equation*}
where $r_\pm=0,1$. Let $\ov{\lambda}=(\ov{\lambda}_i)_{i\geq 1}\in\cP$ be such that $\ov{\lambda}'_1=n-\lambda'_1$ and $\ov{\lambda}'_i=\lambda'_i$ for $i\geq 2$. 
Let $M_+=\lambda'_1$, $M_-=\ov{\lambda}'_1=n-\lambda'_1$, and $L=M_\pm+q_\pm$.
Note that $2L+r_\pm=n$. Then we put
\begin{equation*}\label{product form of T}
\widehat{\bf T}^{\mf d}(\la,n)=
\begin{cases}
{\bf T}^{\mf{d}}(\lambda_{1})\times\cdots \times {\bf T}^{\mf{d}}(\lambda_{M_+})\times {\bf T}^{\mf{d}}(0)^{\times q_+} \times \left({\bf T}^{\rm sp +} \right)^{r_+}, & \text{if $n-2\lambda'_1\geq 0$},\\
{\bf T}^{\mf{d}}(\ov{\lambda}_{1})\times\cdots \times {\bf T}^{\mf{d}}(\ov{\lambda}_{M_-})\times  \ov{\bf T}^{\mf{d}}(0)^{\times q_-}\times \left({\bf T}^{\rm sp -} \right)^{r_-}, & \text{if $n-2\lambda'_1\leq0$}.
\end{cases}
\end{equation*}

\begin{df}\label{def:psst}{\rm For $\la\in {\mc P}({\rm G}_n)$,
we define
\begin{equation*}
{\bf T}^{\mf g}(\la,n)=\{\,{\bf T}=(T_i)_{i\geq 1}\in \widehat{\bf T}^{\mf g}(\la,n)\,|\,\text{ $T_i\prec T_{i+1}$ for $i\geq 1$}\,\}.
\end{equation*}
}
\end{df}

We regard $\widehat{\bf T}^{\mf g}(\la,n)$ as a $\g_\infty$-crystal by identifying ${\bf T}=(T_1,T_2,\ldots)\in \widehat{\bf T}^{\mf g}(\la,n)$ with $\cdots \otimes T_2\otimes T_1$, and regard ${\bf T}^{\mf g}(\la,n)$ as its subcrystal. Let ${\bf H}_\la$ be the element $(T_1,T_2,\ldots)$ in ${\bf T}^{\mf g}(\la,n)$ such that each $T_k$ is a highest weight element for all $k\geq 1$, or $\te_iT_k={\bf 0}$ for all $i\in I$. 
Then we have the following by \cite[Theorem 7.4]{K13} for ${\mf g}={\mf b}, {\mf c}$ and \cite[Theorem 4.4]{K14} for ${\mf g}={\mf d}$.

\begin{thm}
For $\la\in {\mc P}({\rm G}_n)$, we have
\begin{equation*}
{\bf T}^{\mf g}(\la,n) \cong {\bf B}({\mf g}_\infty,\La^{\mf g}(\la)),
\end{equation*}
with highest weight element ${\bf H}_\la$ of weight $\Lambda^{\mf g}(\la)$.
\end{thm}

\begin{rem}\label{rem:case of osp}
{\rm
Following \cite{Je},  one can see that ${\bf T}^{\mf b}(\la,n)$ for $\la\in \mc{P}({\rm G}_n)$ (${\rm G}={\rm Pin}$ or ${\rm Spin}$)  coincides with the crystal graph associated to an integrable highest weight module over $\mf b^\bullet_\infty$ with highest weight $\La^{\mf b}(\la)$. So the tensor product and branching multiplicities for $\mf b_\infty$ are equal to those for $\mf b^\bullet_\infty$. 
Recall that $\mf b^\bullet_k\cong \mf{osp}_{1|2k}$ for $k\in\N$.
}
\end{rem}

\begin{rem}\label{rem:reduction to k-1}
{\rm
For $k\geq 2$ and $\la\in\mc{P}({\rm G}_n)$ with $\la_1\leq k$, let ${\bf T}^{\mf g}_k(\la,n)$ be the subset of ${\bf T}^{\mf g}(\la,n)$ with entries from $\{\,1,\ldots,k\,\}$. Then ${\bf T}^{\mf g}_k(\la,n)$ is a $\mf g_k$-crystal associated to a finite-dimensional irreducible $\mf g_k$-module with highest weight $\La^{\mf g}(\la)$ for $\mf g=\mf b, \mf b^\bullet, \mf c, \mf d$ with respect to $\te_i$ and $\tf_i$ for $i\in\{0,\ldots,k-1\}$.
}
\end{rem}

\section{Stability in tensor product and branching decomposition}

\subsection{Separation Lemma}
We put
\begin{equation}\label{partition of type x}
\begin{split}
\cP_{\mf b}&=\cP,\\
\cP_{\mf c}&=\{\,2\la=(2\la_i)_{i\geq 1}\,|\,\la\in \cP\,\},\\
\cP_{\mf d}&=\{\,(2\la)'\,|\,\la\in \cP\,\}.
\end{split}
\end{equation}

Consider a ${\mf g}_\infty$-crystal ${\bf T}^{\mf g}(\la,n)$ for $\la\in \mc{P}({\rm G}_n)$. Let ${\bf T}=(T_1,\ldots,T_r)\in {\bf T}^{\mf g}(\la,n)$ given. We put $w({\bf T})=w(T_r)\ldots  w(T_1)$, and define $L({\bf T})$ to be the maximal length of a weakly decreasing subword of $w({\bf T})$.

\begin{lem}\label{lem:aux-2}
For ${\bf T} \in {\bf T}^{\mf g}(\la,n)$, we have $L({\bf T})\geq L({\bf H}_\la)=\ell(\la)$.
\end{lem}
\pf It is clear that $L({\bf H}_\la)=\ell(\la)$. For arbitrary ${\bf T}$, it follows directly from Lemma \ref{lem:aux-1}(1) that $L({\bf T})\geq \ell(\la)$.
\qed

We put
\begin{equation}\label{eq:body and tail}
{\bf T}^{\tt tail}=(T^{\tt tail}_1,\ldots, T^{\tt tail}_r), \quad 
{\bf T}^{\tt body}=(T^{\tt body}_1,\ldots, T^{\tt body}_r).
\end{equation}
For ${\mf g}={\mf d}$ and $T_i\in \ov{\bf T}^{\mf d}(0)$ or ${\bf T}^{\rm sp -}$, we assume that $T_i^{\tt tail}$ is the subtableau  of $T_i$ consisting of the boxes at the bottom, and $T_i^{\tt body}$ is the subtableau consisting of the other part of $T_i$.

We assume that $T_i^{\tt body}$ and $T_i^{\tt tail}$ are separated by a common horizontal line. Then we may regard ${\bf T}^{\tt tail}$ as a tableau of shape $\la'$, where $T_i^{\tt tail}$ is the $i$th column from the left for $1\leq i\leq r$, and regard ${\bf T}^{\tt body}$ as a tableau of shape $\mu^\pi$ for some $\mu\in\cP$, when ${\mf r}_{T_i}=0$ for $1\leq i\leq r$.
Here $\mu^\pi$ denotes a skew Young diagram obtained by $180^\circ$-rotation of $\mu$.

\begin{ex}{\rm
Let ${\bf T}=(T_1,T_2)\in {\bf T}^{\mf b}(2,1)$ be given by

$$ T_1 = \  \resizebox{.055\hsize}{!}
{\def\lr#1{\multicolumn{1}{|@{\hspace{.75ex}}c@{\hspace{.75ex}}|}{\raisebox{-.04ex}{$#1$}}}\raisebox{-.6ex}
{$\begin{array}{cc}
\\
\cline{2-2}
&\lr{\color{red}4}\\
\cline{1-1}\cline{2-2}
\lr{\color{red}5} & \lr{\color{red}8}\\
\cline{1-1}\cline{2-2}
\lr{\color{blue}6} \\
\cline{1-1}
\lr{\color{blue}7} \\
\cline{1-1}
\end{array}$}} 
  \ \ \ \ \    T_2 = \  \resizebox{.055\hsize}{!}
{\def\lr#1{\multicolumn{1}{|@{\hspace{.75ex}}c@{\hspace{.75ex}}|}{\raisebox{-.04ex}{$#1$}}}\raisebox{-.6ex}
{$\begin{array}{cc}
\cline{1-1}\cline{2-2}
\lr{\color{red}3} &\lr{\color{red}5}\\
\cline{1-1}\cline{2-2}
\lr{\color{red}4}&\lr{\color{red}8}\\
\cline{1-1}\cline{2-2}
\lr{\color{red}7} & \lr{\color{red}9}\\
\cline{1-1}\cline{2-2}
\lr{\color{blue}8} \\
\cline{1-1}\\ 
\end{array}$}}$$
where the entries in $T_i^{\tt body}$ and $T_i^{\tt tail}$ ($i=1,2$) are in red and blue, respectively. 
Then we have
$$ 
 {\bf T}^{\texttt{tail}} = \quad  \resizebox{.055\hsize}{!}
{\def\lr#1{\multicolumn{1}{|@{\hspace{.75ex}}c@{\hspace{.75ex}}|}{\raisebox{-.04ex}{$#1$}}}\raisebox{-.6ex}
{$\begin{array}{cc}
\\ 
\cline{1-1}\cline{2-2}
\lr{\color{blue}6} & \lr{\color{blue}8}\\
\cline{1-1}\cline{2-2}
\lr{\color{blue}7} \\
\cline{1-1}
\end{array}$}} \ \ \ \ \ \ \  
{\bf T}^{\texttt{body}} = \  \resizebox{.11\hsize}{!}
{\def\lr#1{\multicolumn{1}{|@{\hspace{.75ex}}c@{\hspace{.75ex}}|}{\raisebox{-.04ex}{$#1$}}}\raisebox{-.6ex}
{$\begin{array}{cccc}
\cline{3-3}\cline{4-4}
& & \lr{\color{red}3} &\lr{\color{red}5}\\
\cline{2-2}\cline{3-3}\cline{4-4}
&\lr{\color{red}4} & \lr{\color{red}4}&\lr{\color{red}8}\\
\cline{1-1}\cline{2-2}\cline{3-3}\cline{4-4}
\lr{\color{red}5} & \lr{\color{red}8} & \lr{\color{red}7} & \lr{\color{red}9}\\
\cline{1-1}\cline{2-2}\cline{3-3}\cline{4-4}\\ 
\end{array}$}} $$\vskip 3mm
\noindent Note that ${\bf T}^{\tt body}$ is not semistandard.

}
\end{ex}

The following lemma will play a crucial role in this paper.

\begin{lem}\label{lem:aux-3}
If $L({\bf T})\leq \frac{n}{2}$, then we have
\begin{itemize}
\item[(1)] ${\bf T}^{\tt tail}\in SST(\la')$, and ${\bf T}^{\tt body}\in SST(\mu^\pi)$ for some $\mu\in\cP_{\mf g}$,

\item[(2)] ${\bf T}\equiv_A {\bf T}^{\tt body}\otimes {\bf T}^{\tt tail}$, where $\equiv_A$ denotes the equivalence as elements of $\mf l_\infty$-crystals.
\end{itemize}
\end{lem}
\pf \textsc{Case 1}. Suppose that ${\mf g}={\mf c}$. Let ${\bf T}=(T_1,\ldots, T_{n/2})\in {\bf T}^{\mf c}(\la,n)$ with $L({\bf T})\leq n/2$.
We may regard ${\bf T}^{\tt tail}$ as a tableau of shape $\la'$ and ${\bf T}^{\tt body}$ as a tableau of shape $\mu^\pi$ for some $\mu\in \cP_{\mf c}$.
We first have by Lemma \ref{lem:aux-1}(1) that ${\bf T}^{\tt tail}$ is semistandard and hence ${\bf T}^{\tt tail}\in SST(\la')$.

Next, we claim that ${\bf T}^{\tt body}$ is semistandard. Since $T_i^{\tt body}$ is a semistandard tableau of shape $\la(0,0,c_i)$ for some $c_i\in \Z_+$ with $(T_i^{\tt body})^{\tt R}=T_i^{\tt R}$ and $(T_i^{\tt body})^{\tt L}(k)=T_i^{\tt L}(\la_i+k)$ for $1\leq i\leq n/2$ and $k\geq 1$, it is enough to show that  $(T_i^{\tt body})^{\tt R}(k)\leq (T_{i+1}^{\tt body})^{\tt L}(k)$ for all $i$ and $k$.

Suppose that $(T_i^{\tt body})^{\tt R}(k_0)> (T_{i+1}^{\tt body})^{\tt L}(k_0)$ for some $i$ and  $k_0$.  Since $(T_i^{\tt body})^{\tt R}(k_0)> (T_{i+1}^{\tt body})^{\tt L}(k_0)$ and $T_i\prec T_{i+1}$, at least one of the boxes in $(T_{i+1}^{\tt body})^{\tt L}$ moves to the right to form an entry in ${}^{\tt R}T_{i+1}$. This implies that $T_{i+1}^{\tt tail}$ is not empty, and  $T_{i+1}^{\tt tail}[k_1]\leq (T_{i+1}^{\tt body})^{\tt R}(1)$ for some $k_1\geq 1$.
By Lemma \ref{lem:aux-1} we have a weakly increasing sequence
\begin{equation*}
\begin{split}
T_{1}^{\tt tail}[k_1]&\leq T_{2}^{\tt tail}[k_1]\leq  \ldots\leq T_{i}^{\tt tail}[k_1]\leq T_{i+1}^{\tt tail}[k_1]
\\
&\leq T_{i+1}^{\tt R}(1)\leq T_{i+2}^{\tt R}(1)\leq \ldots \leq T_{n/2}^{\tt R}(1).
\end{split}
\end{equation*}
which yields a weakly decreasing subword of $w({\bf T})$ of length $n/2+1$. This is a contradiction, which proves our claim. Therefore, ${\bf T}^{\tt body}\in SST(\mu^\pi)$.

Finally, we claim that ${\bf T}\equiv_A {\bf T}^{\tt body}\otimes {\bf T}^{\tt tail}$. By the semistandardness of ${\bf T}^{\tt body}$ and Lemma \ref{lem:aux-1}(2), we have for $2\leq i\leq n/2$
\begin{equation}\label{eq:aux-1}
\ldots \leq(T_{i-1}^{\tt body})^{\tt L}(1)\leq (T_{i-1}^{\tt body})^{\tt R}(1)\leq (T_i^{\tt body})^{\tt L}(1)<T_i^{\tt tail}[1].
\end{equation}
Since $\equiv_A$ coincides with the usual Knuth equivalence,
one can check by \eqref{eq:aux-1} and considering the column insertion of tableaux that
\begin{equation}\label{eq:aux-2}
\begin{split}
T_i^{\tt body}\otimes T_i^{\tt tail}\otimes &\left(T_{i-1}^{\tt body}\otimes \cdots \otimes T_{1}^{\tt body}\right) \\ &\equiv_A
T_i^{\tt body}\otimes \left(T_{i-1}^{\tt body}\otimes \cdots \otimes T_{1}^{\tt body}\right)\otimes T_i^{\tt tail}.
\end{split}
\end{equation}
Applying \eqref{eq:aux-2} successively, we have
\begin{equation*}
\begin{split}
{\bf T} & = T_{n/2}\otimes \cdots \otimes T_1 \\
& \equiv_A \left(T_{n/2}^{\tt body}\otimes T_{n/2}^{\tt tail}\right) \otimes \cdots \otimes  \left(T_{1}^{\tt body}\otimes T_{1}^{\tt tail}\right)\\
& \equiv_A \left(T_{n/2}^{\tt body}\otimes \ldots\otimes T_{1}^{\tt body}\right) \otimes \left( T_{n/2}^{\tt tail}\otimes \cdots \otimes T_{1}^{\tt tail}\right)\equiv_A {\bf T}^{\tt body} \otimes {\bf T}^{\tt tail}.
\end{split}
\end{equation*}\vskip 2mm

\textsc{Case 2}. Suppose that ${\mf g}={\mf b}$. The proof is almost identical to \textsc{Case 1}. We leave the details to the readers.
\vskip 2mm

\textsc{Case 3}. Suppose that ${\mf g}={\mf d}$. Let ${\bf T}=(T_1,\ldots, T_{r})\in {\bf T}^{\mf d}(\la,n)$ with $L({\bf T})\leq n/2$. By Lemma \ref{lem:aux-2}, we have $\ell(\lambda)\leq n/2$, and hence $n-2\lambda'_1\geq 0$.

We first claim that ${\mf r}_{T_i}=0$ for all $1\leq i\leq r$. Suppose that ${\mf r}_{T_i}=1$ for some $1\leq i\leq \lambda'_1$, and let $j$ be the smallest one such that ${\mf r}_{T_j}=1$. By Lemma \ref{lem:aux-1}, we have a sequence
\begin{equation*}
T_1^{\tt tail}[1]\leq  \ldots \leq T_{j}^{\tt tail}[1] \leq T_{j}^{\tt R}(1)\leq T_{j+1}^{\tt R}(1)\leq \ldots \leq T_{r}^{\tt R}(1).
\end{equation*}
Note that $T_{j+1}^{\tt R}(1), \ldots, T_{r}^{\tt R}(1)$ are non-empty since $T_i \prec T_{i+1}$ for $j\leq i\leq r-1$.
Reading from right to left, we get a subword of $w({\bf T})$ with length $r+1$. Since $2\lambda'_1+2q_++r_+=n$ and $r=\lambda'_1+q_++r_+$, we have $r+1>n/2$. This is a contradiction, which proves our claim.

Now by Lemma \ref{lem:aux-1}(1), we have ${\bf T}^{\tt tail}\in SST(\la')$.
By the same arguments as in \textsc{Case 1} we conclude that ${\bf T}^{\tt body}\in SST(\mu^\pi)$ and ${\bf T}\equiv_A {\bf T}^{\tt body}\otimes {\bf T}^{\tt tail}$.
\qed

\subsection{Stable tensor product and branching rules}
Let $\mu\in {\mc P}({\rm G}_m)$ and $\nu\in {\mc P}({\rm G}_{n})$ be given. For $\la\in \mc{P}({\rm G}_{m+n})$, let
\begin{equation}\label{eq:LR-1}
\begin{split}
\texttt{\large LR}^{\la}_{\mu\nu }({\mf g})&=
\left\{\,{\bf T}\in {\bf T}^{\mf g}(\nu,n)\,\Bigg \vert \,
\begin{array}{l}
(1)\ {\rm wt}({\bf T})+{\rm wt}({\bf H}_\mu)={\rm wt}({\bf H}_\la)\\
(2)\ \varepsilon_i({\bf T})\leq \langle \Lambda^{\mf g}(\mu), \alpha^\vee_i \rangle \ \ \text{for $i\in I$}
\end{array}
\,\right\}, \\
c^{\la}_{\mu\nu }({\mf g})
&= \left| \texttt{\large LR}^{\la}_{\mu\nu }({\mf g}) \right|.
\end{split}
\end{equation}
Then we have the following formula of tensor product multiplicity for integrable highest weight ${\mf g}_\infty$-modules.
\begin{prop}\label{prop:tensor product}
For $\mu\in {\mc P}({\rm G}_m)$, $\nu\in {\mc P}({\rm G}_n)$, and $\la\in \mc{P}({\rm G}_{m+n})$, we have
\begin{equation*}
c^{\la}_{\mu\nu }({\mf g})=
\left[
L({\mf g}_\infty,\Lambda^{\mf g}(\mu))\otimes L({\mf g}_\infty,\Lambda^{\mf g}(\nu)):L({\mf g}_\infty,\Lambda^{\mf g}(\la))
\right].
\end{equation*}
\end{prop}
\pf By tensor product rule of crystals, we have $\texttt{\large LR}^{\la}_{\mu\nu }({\mf g})=\{\,{\bf T}\in {\bf T}^{\mf g}(\nu,n)\,|\,{\bf H}_\mu \otimes {\bf T}\equiv {\bf H}_\la \,\}$, and hence $c^{\la}_{\mu\nu }({\mf g})$ is the multiplicity of ${\bf B}({\mf g}_\infty,\La^{\mf g}(\la))$ in ${\bf B}({\mf g}_\infty,\La^{\mf g}(\mu))\otimes {\bf B}({\mf g}_\infty,\La^{\mf g}(\nu))$. \qed\vskip 2mm

Now we have the following stable tensor product rule, which is the first main result in this paper.

\begin{thm}\label{thm:stable tensor}
Let $\mu\in {\mc P}({\rm G}_m)$, $\nu\in {\mc P}({\rm G}_n)$, and $\la\in \mc{P}({\rm G}_{m+n})$ be given. If $\ell(\la)\leq \frac{1}{2}\min\{m,n\}$, then the map ${\bf T}\mapsto ({\bf T}^{\tt body}, {\bf T}^{\tt tail})$ gives a bijection
\begin{equation*}
{\tt \large LR}^{\la}_{\mu\nu }({\mf g}) \longrightarrow \bigsqcup_{\gamma\in\cP}\bigsqcup_{\delta\in \cP_{\mf g}} {\tt LR}^{\gamma}_{\mu' \delta}\times {\tt LR}^{\la'}_{\gamma\nu'}.
\end{equation*}
\end{thm}
\pf Let ${\bf T}\in \texttt{\large LR}^{\la}_{\mu\nu }({\mf g})$ be given. By definition, ${\bf H}_\mu \otimes {\bf T}\equiv {\bf H}_\la$ as a ${\mf g}_\infty$-crystal element, and hence as  an ${\mf l}_\infty$-crystal element. Recall that $\equiv_{A}$ denotes the equivalence of ${\mf l}_\infty$-crystal elements. Note that ${\bf H}_\la \equiv_A H_{\la'}$.
Since $\ell(\la)\leq \frac{1}{2}\min\{m,n\}$, we have  $L({\bf T})\leq L({\bf H}_\mu \otimes {\bf T})=\ell(\la)\leq \frac{1}{2}\min\{m,n\}\leq n/2$. Here we understand $L({\bf H}_\mu \otimes {\bf T})$ as the maximal length of a weakly decreasing subword of $w({\bf H}_\mu)w({\bf T})$.
By Lemma \ref{lem:aux-3}, we have
\begin{itemize}
\item[$\cdot$] ${\bf T}^{\tt body}\in SST(\delta^\pi)$ for some $\delta\in \cP_{\mf g}$,

\item[$\cdot$] ${\bf T}^{\tt tail}\in SST(\nu')$,

\item[$\cdot$] ${\bf T}\equiv_A {\bf T}^{\tt body}\otimes {\bf T}^{\tt tail}$.

\end{itemize}
It is clear that the correspondence ${\bf T}\mapsto ({\bf T}^{\tt body}, {\bf T}^{\tt tail})$ is injective.

Since ${\bf H}_\mu\otimes {\bf T}^{\tt body}\equiv_A H_{\mu'}\otimes {\bf T}^{\tt body}\equiv_A H_\gamma$ for some $\gamma\in \cP$  and there exists a unique tableau $U\in SST(\delta)$ such that $U\equiv_A {\bf T}^{\tt body}$, we may regard ${\bf T}^{\tt body}$ as an element in ${\tt LR}^\gamma_{\mu' \delta}$ (see Section \ref{subsec:tableaux}).
Then we have $H_\gamma\otimes {\bf T}^{\tt tail}\equiv_A {\bf H}_\lambda\equiv_A H_{\lambda'}$, which implies that ${\bf T}^{\tt tail}\in {\tt LR}^{\la'}_{\gamma \nu'}$.

Conversely, suppose that a pair $(U,V)\in {\tt LR}^{\gamma}_{\mu' \delta}\times {\tt LR}^{\la'}_{\gamma\nu'}$ is given for $\gamma\in \cP$ and $\delta\in \cP_{\mf g}$. Here we assume that $U\in SST(\delta^\pi)$ such that $H_{\mu'}\otimes U\equiv_A H_\gamma$, and $V\in SST(\nu')$ such that $H_{\gamma}\otimes U\equiv_A H_{\lambda'}$. We also assume that ${\mf g}={\mf c}$ since the arguments for ${\mf g}={\mf b}, {\mf d}$ are very similar.

Note that $\ell(\delta'), \ell(\nu)\leq \ell(\la)\leq \frac{1}{2}\min\{m,n\}$.
Let $U^{(n-i+1)}$ be the $i$th column of $U$ from the right for $1\leq i\leq n$, where $U^{(i)}$ is empty for $1\leq i\leq n/2$. Let $V^{(i)}$ be the $i$th column of $V$ from the left for $i\geq 1$. 
We claim that for all $1\leq i\leq n/2$
\begin{equation}\label{eq:UV admissible}
V^{(i)}[1]> U^{(2i)}(1).
\end{equation}
Otherwise, we have a weakly decreasing subword of $w(U)w(V)$ given by
\begin{equation*}
U^{(n)}(1)\geq U^{(n-1)}(1)\geq \ldots\geq U^{(2i)}(1)\geq V^{(i)}[1]\geq V^{(i-1)}[1]\geq \ldots \geq V^{(1)}[1],
\end{equation*}
for some $i$, which is of length $(n-2i+1)+i =n-i+1 > n/2$. Since $H_{\mu'}\otimes U\otimes V \equiv_A H_{\la'}$, this implies that $\ell(\la)> n/2$, which is a contradiction.

For $1\leq i\leq n/2$, let $T_i$ be the tableau of shape $\la(\nu_i,0,c)$ for some $c\in\Z_+$ such that $(T_i^{\tt body})^{\tt L}=U^{(2i-1)}$, $(T_i^{\tt body})^{\tt R}=U^{(2i)}$, and $T_i^{\tt tail}=V^{(i)}$. By \eqref{eq:UV admissible}, it is straightforward to see that $T_i\in {\bf T}^{\mf g}(\nu_i)$ and $T_i\prec T_{i+1}$. Therefore, ${\bf T}=(T_1,\ldots, T_{n/2})\in {\bf T}^{\mf g}(\nu,n)$ with ${\bf T}^{\tt body}=U$ and ${\bf T}^{\tt tail}=V$. 
By Lemma \ref{lem:aux-3}, we have ${\bf H}_\mu\otimes {\bf T}\equiv_A H_{\mu'}\otimes U\otimes V \equiv_A H_{\lambda'}$, and hence $\te_i({\bf H}_\mu\otimes {\bf T})=0$ for $i\in I\setminus \{0\}$.

We see from the definition of $\td{e}_0$ and $\td{f}_0$ that $\varphi_0({\bf H}_\mu)= \frac{m}{2}-\ell(\mu)$, and $\varepsilon_0({\bf T})\leq z/2$, where $z$ is the number of $1$'s in $U$. Since $H_{\mu'}\otimes U \otimes V \equiv_A H_{\lambda'}$, we have $z  + \ell(\mu) \leq \ell(\lambda)$, and then
\begin{equation*}
\varepsilon_0({\bf T}) + \ell(\mu) \leq z+\ell(\mu) \leq \ell(\lambda) \leq \frac{m}{2},
\end{equation*}
which implies that $\varepsilon_0({\bf T}) \leq \frac{m}{2}-\ell(\mu)=\varphi_0({\bf H}_{\mu})$, and $\td{e}_0({\bf H}_\mu\otimes {\bf T})=0$. Hence ${\bf T}\in {\tt \large LR}^{\la}_{\mu\nu }({\mf g})$.

Therefore, it follows that the map ${\bf T}\mapsto ({\bf T}^{\tt body}, {\bf T}^{\tt tail})$ gives a bijection.
\qed

\begin{cor}\label{cor:stable tensor}
Let $\mu\in {\mc P}({\rm G}_m)$, $\nu\in {\mc P}({\rm G}_n)$, and $\la\in \mc{P}({\rm G}_{m+n})$ be given.
If $\ell(\la)\leq \frac{1}{2}\min\{m,n\}$, then
\begin{equation*}
c^{\la}_{\mu\nu }({\mf g})=\sum_{\gamma\in \cP}\sum_{\delta\in \cP_{\mf g}}c^{\gamma}_{\mu \delta'}c^{\la}_{\gamma\nu}.
\end{equation*}
\end{cor}
\pf It follows directly from Theorem \ref{thm:stable tensor} and the invariance of Littlewood-Richardson coefficients under conjugate of partitions.
\qed\vskip 2mm

Let $\mu\in {\mc P}({\rm G}_n)$ be given. For $\la\in\cP$, let
\begin{equation}\label{eq:LR-2}
\begin{split}
\texttt{\large LR}^{\mu}_{\la}({\mf g})&=
\left\{\,{\bf T}\in {\bf T}^{\mf g}(\mu,n)\,\Bigg \vert \,
\begin{array}{l}
(1)\ {\rm wt}({\bf T})= \frac{n}{\epsilon}\La^{\mf g}_0+\sum_{i\geq 1}\la'_i\epsilon_i \\
(2)\ \varepsilon_i({\bf T})=0 \ \ \text{for $i\in I\setminus\{0\}$}
\end{array}
\,\right\},\\
c^\mu_{\la}({\mf g})&=\left| \texttt{\large LR}^\mu_{\la}({\mf g}) \right|.
\end{split}
\end{equation}

\begin{prop}\label{prop:branching to levi}
For $\mu\in {\mc P}({\rm G}_n)$ and $\la\in \cP$, we have
\begin{equation*}
c^\mu_{\la}({\mf g})=
\left[
L({\mf g}_\infty,\Lambda^{\mf g}(\mu)):L({\mf l}_\infty,\la')
\right],
\end{equation*}
where $L({\mf l}_\infty,\la')$ is the irreducible ${\mf l}_\infty$-module with highest weight $\sum_{i\geq 1}\la'_i\epsilon_i$.
\end{prop}
\pf By definition, the connected component of ${\bf T}\in \texttt{\large LR}^\mu_{\la}({\mf g})$ in ${\bf T}^\g(\mu,n)$ with respect to $\te_i$ and $\tf_i$ for $i\in I\setminus\{0\}$ is isomorphic to $SST(\la')$ as an ${\mf l}_\infty$-crystal. Hence $c^\mu_{\la}({\mf g})$ is the multiplicities of $SST(\la')$ or $\B({\mf l}_\infty, \la')$ in ${\bf B}({\mf g}_\infty,\La^{\mf g}(\la))$.
\qed\vskip 2mm

Then we have  following stable branching rule, which is the second main result in this paper.

\begin{thm}\label{thm:stable branching}
Let $\mu\in {\mc P}({\rm G}_n)$ and $\la\in \cP$ be given. If $\ell(\la) \leq \frac{n}{2}$, then the map ${\bf T}\mapsto {\bf T}^{\tt tail}$ gives a bijection 
\begin{equation*}
{\tt \large LR}^\mu_{\la}({\mf g}) \longrightarrow  \bigsqcup_{\delta\in \cP_{\mf g}} {\tt LR}^{\la'}_{\delta\mu'}.
\end{equation*}
\end{thm}
\pf Let ${\bf T}\in \texttt{\large LR}^\mu_{\la}({\mf g})$ be given. By definition, ${\bf T}\equiv_A H_{\la'}$.
Since $L({\bf T})=\ell(\la)\leq \frac{n}{2}$, we have  by Lemma \ref{lem:aux-3}
\begin{itemize}
\item[$\cdot$] ${\bf T}^{\tt body}\in SST(\delta^\pi)$ for some $\delta\in \cP_{\mf g}$,

\item[$\cdot$] ${\bf T}^{\tt tail}\in SST(\mu')$,

\item[$\cdot$] ${\bf T}\equiv_A {\bf T}^{\tt body}\otimes {\bf T}^{\tt tail}$.

\end{itemize}
Moreover, since ${\bf T}\equiv_A {\bf T}^{\tt body}\otimes {\bf T}^{\tt tail}\equiv_A H_{\la'}$, we have ${\bf T}^{\tt body}\equiv_A H_{\delta}$ and we may regard ${\bf T}^{\tt tail}$ as an element in ${\tt LR}^{\la'}_{\delta \mu'}$.
Hence ${\bf T}\mapsto {\bf T}^{\tt tail}$ is an injective map from ${\tt \large LR}^{\mu}_{\la}({\mf g})$ to $\bigsqcup_{\delta\in \cP_{\mf g}} {\tt LR}^{\la'}_{\delta\mu'}$.

Conversely, suppose that $V\in {\tt LR}^{\la'}_{\delta\mu'}$ is given for $\delta\in \cP_{\mf g}$. Let us assume that ${\mf g}={\mf c}$ since the arguments for ${\mf g}={\mf b}, {\mf d}$ are similar.
Let $U$ be the unique tableau in $SST(\delta^\pi)$ such that $U\equiv_A H_{\delta}$. 
We assume that $V\in SST(\mu')$ such that $U\otimes V\equiv_A H_{\la'}$.
Let $U^{(n-i+1)}$ be the $i$th column of $U$ from the right, and let $V^{(i)}$ be the $i$th column of $V$ from the left for $1\leq i\leq n$. Then $U^{(i)}$ is empty for $1\leq i\leq n/2$, and $V^{(i)}$ is empty for $n/2<i \leq n$ since $\ell(\delta'), \ell(\mu)\leq \ell(\la)\leq n/2$.
By the same argument as in Theorem \ref{thm:stable tensor}, we can show that
$V^{(i)}[1]> U^{(2i)}(1)$ for all $1\leq i\leq n/2$, and hence there exists a unique ${\bf T}\in {\bf T}^{\mf g}(\mu,n)$ such that ${\bf T}^{\tt body}=U$ and ${\bf T}^{\tt tail}=V$. Since ${\bf T}\equiv_A U\otimes V\equiv_A H_{\la'}$, we have ${\bf T}\in {\tt \large LR}^\mu_{\la}({\mf g})$.

Therefore ${\bf T}\mapsto {\bf T}^{\tt tail}$ gives a bijection.
\qed

\begin{cor}\label{cor:stable branching}
Let $\mu\in {\mc P}({\rm G}_n)$ and $\la\in \cP$.
If $\ell(\la)\leq \frac{n}{2}$, then
\begin{equation*}
c^\mu_{\la}({\mf g})= \sum_{\delta\in \cP_{\mf g}} c^{\la}_{\delta'\mu}.
\end{equation*}
\end{cor}

\begin{rem}\label{rem:reduction to k}
{\rm 
For $k\geq 2$ and $\la\in\mc{P}({\rm G}_n)$ with $\la_1\leq k$, let ${\bf T}^{\mf g}_k(\la,n)$ be as in Remark \ref{rem:reduction to k-1}.
Since ${\bf T}^{\mf g}_k(\la,n)$ for $\la\in \bigsqcup_{n}\mc{P}({\rm G}_n)$ with $\la_1\leq k$ give crystals for all finite-dimensional irreducible $\mf g_k$-modules, we have combinatorial extensions of the stable branching rules for the pairs ${\rm G}_\ell \times {\rm G}_\ell  \supset {\rm G}_\ell$ and ${\rm G}_\ell \supset {\rm GL}_{[\ell/2]}$, where ${\rm G}_\ell={\rm Sp}_\ell$, ${\rm Spin}_\ell$, ${\rm Osp}_{1|\ell}$ with $k=[\ell/2]$
by Theorems \ref{thm:stable tensor} and \ref{thm:stable branching} (see also Remark \ref{rem:case of osp}). 
}
\end{rem}

\section{Branching rules for classical groups}


\subsection{Branching rule for ${\rm G}_{m+n}\supset {\rm G}_m\times {\rm G}_n$}
Suppose that ${\rm G}={\rm Sp}$ or ${\rm O}$, that is, ${\mf g}={\mf c}$ or ${\mf d}$.
Note that for $\la\in \mc{P}({\rm G}_{m+n})$, $V_{{\rm G}_{m+n}}^{\la}$ can be viewed as a ${\rm G}_m\times {\rm G}_n$-module since there is a natural embedding of ${\rm G}_m\times {\rm G}_n$ into ${\rm G}_{m+n}$. For $\mu\in {\mc P}({\rm G}_m)$ and $\nu\in {\mc P}({\rm G}_n)$, let $[V_{{\rm G}_{m+n}}^{\la}: V_{{\rm G}_m}^\mu \otimes V_{{\rm G}_n}^{\nu}]$ denote the multiplicity of $ V_{{\rm G}_m}^\mu \otimes V_{{\rm G}_n}^{\nu}$ in $V_{{\rm G}_{m+n}}^{\la}$.

\begin{thm}\label{thm:stable tensor-dual}
For $\la\in \mc{P}({\rm G}_{m+n})$, $\mu\in {\mc P}({\rm G}_m)$ and $\nu\in {\mc P}({\rm G}_n)$, we have
\begin{equation*}
\left[V_{{\rm G}_{m+n}}^{\la}: V_{{\rm G}_m}^\mu \otimes V_{{\rm G}_n}^{\nu}\right]=c^\la_{\mu\nu}({\mf g}),
\end{equation*}
where $c^\la_{\mu\nu}({\mf g})$ is given in \eqref{eq:LR-1}.
Furthermore, if $\ell(\la)\leq \frac{1}{2}\min\{m,n\}$, then
\begin{equation*}
\left[V_{{\rm G}_{m+n}}^{\la}: V_{{\rm G}_m}^\mu \otimes V_{{\rm G}_n}^{\nu}\right]=\sum_{\gamma\in \cP}\sum_{\delta\in \cP_{\mf g}}c^{\gamma}_{\mu \delta'}c^{\la}_{\gamma\nu}.
\end{equation*}
\end{thm}
\pf By Propositions \ref{duality} and \ref{prop:tensor product}, we have
{\allowdisplaybreaks
\begin{equation*}
\begin{split}
\mc{F}^{\frac{m}{2}}\otimes \mc{F}^{\frac{n}{2}}&\cong
\left( \bigoplus_{\mu\in {\mc P}({\rm G}_m)}
L(\mf{g}_\infty,\La^{\mf g}(\mu))\otimes V_{{\rm G}_m}^\mu \right)\otimes
\left( \bigoplus_{\nu\in {\mc P}({\rm G}_n)}
L(\mf{g}_\infty,\La^{\mf g}(\nu))\otimes V_{{\rm G}_n}^\nu \right)\\
&\cong
\bigoplus_{\mu\in {\mc P}({\rm G}_m)}
\bigoplus_{\nu\in {\mc P}({\rm G}_n)}
\left(L(\mf{g}_\infty,\La^{\mf g}(\mu))\otimes L(\mf{g}_\infty,\La^{\mf g}(\nu))\right)\otimes
\left(  V_{{\rm G}_m}^\mu\otimes V_{{\rm G}_n}^\nu \right)\\
&\cong
\bigoplus_{\mu\in {\mc P}({\rm G}_m)}
\bigoplus_{\nu\in {\mc P}({\rm G}_n)}
\left(\bigoplus_{\la\in \mc{P}(G_{m+n})}L(\mf{g}_\infty,\La^{\mf g}(\la))^{\oplus c^\la_{\mu\nu}({\mf g})} \right)\otimes
\left(  V_{{\rm G}_m}^\mu\otimes V_{{\rm G}_n}^\nu \right)\\
&\cong \bigoplus_{\la\in {\mc P}({\rm G}_{m+n})}
L(\mf{g}_\infty,\La^{\mf g}(\la))\otimes \left( \bigoplus_{\mu\in {\mc P}({\rm G}_m)}
\bigoplus_{\nu\in {\mc P}({\rm G}_n)} \left(V_{{\rm G}_m}^\mu\otimes V_{{\rm G}_n}^\nu\right)^{\oplus c^\la_{\mu\nu}({\mf g})} \right).
\end{split}
\end{equation*}}
Since $\mc{F}^{\frac{m}{2}}\otimes \mc{F}^{\frac{n}{2}}\cong \mc{F}^{\frac{m+n}{2}}$, we have by Proposition \ref{duality}
\begin{equation*}
\left[V_{{\rm G}_{m+n}}^{\la}: V_{{\rm G}_m}^\mu \otimes V_{{\rm G}_n}^{\nu}\right]=c^\la_{\mu\nu}({\mf g}).
\end{equation*}
The second equality follows from Corollary \ref{cor:stable tensor}.
\qed\vskip 3mm

\begin{rem}{\rm
By the same argument as in Theorem \ref{thm:stable tensor-dual}, we also have a branching rule for the closed subgroup ${\rm Spin}_m\times {\rm Spin}_n\subset {\rm Pin}_{m+n}$ for $m,n$ odd in the case when $\mf g=\mf b$. 
}
\end{rem}

\subsection{Branching rule for ${\rm GL}_{n}\supset {\rm G}_n$}
Suppose that ${\rm G}={\rm Sp}$ or ${\rm O}$. Note that ${\rm G}_n$ is a subgroup of ${\rm GL}_n$. For $\la\in \cP$ with $\ell(\la)\leq n$, let $V^\la_{{\rm GL}_n}$ be the finite-dimensional irreducible ${\rm GL}_n$-module corresponding to $\la$, and for $\mu\in \mc{P}({\rm G}_n)$, let $\left[ V^\la_{{\rm GL}_n} : V^\mu_{{\rm G}_n} \right]$ denote the multiplicity of $V^\mu_{{\rm G}_n}$ in $V^\la_{{\rm GL}_n}$ as a ${\rm G}_n$-module. Then we have the following, which extends the Littlewood restriction rule \cite{Lw-1,Lw-2}.

\begin{thm}\label{thm:stable branching-dual}
For $\la\in \cP$ with $\ell(\la)\leq n$ and $\mu\in \mc{P}({\rm G}_n)$, we have
\begin{equation*}
\left[ V^\la_{{\rm GL}_n} : V^\mu_{{\rm G}_n} \right]=c^\mu_\la({\mf g}),
\end{equation*}
where $c^\mu_\la({\mf g})$ is given in \eqref{eq:LR-2}.
Furthermore, if $\ell(\la)\leq \frac{n}{2}$, then
\begin{equation*}
\left[ V^\la_{{\rm GL}_n} : V^\mu_{{\rm G}_n} \right] = \sum_{\delta\in \cP_{\mf g}} c^{\la}_{\delta'\mu}.
\end{equation*}
\end{thm}
\pf We first consider an action of $\gl_n=\bigoplus_{1\leq i,j\leq n}\mathbb{C}e_{ij}$ as operators on $\mc{F}^{\frac{n}{2}}$ by left multiplication as follows:

For $n=2\ell$ and ${\rm G}_n={\rm O}_{2\ell}$, put
\begin{equation}\label{eq:gl-1}
e_{ij}=
\begin{cases}
\sum_{r\in \hf+\Z_+}\psi^{+,i}_{-r}\psi^{-,j}_{r}, & \text{if $1\leq i, j\leq \ell$},\\
\sum_{r\in \hf+\Z_+}\psi^{+,i}_{-r}\psi^{+,j-\ell}_{r}, & \text{if $1\leq i\leq \ell$ and $\ell<j\leq n$},\\
\sum_{r\in \hf+\Z_+}\psi^{-,i-\ell}_{-r}\psi^{-,j}_{r}, & \text{if $\ell<i\leq n$ and $1\leq j\leq \ell$},\\
\sum_{r\in \hf+\Z_+}\psi^{-,i-\ell}_{-r}\psi^{+,j-\ell}_{r}, & \text{if $\ell<i,j\leq n$}.
\end{cases}
\end{equation}

For $n=2\ell$ and ${\rm G}_n={\rm Sp}_{2\ell}$, put
\begin{equation}\label{eq:gl-2}
e_{ij}=
\begin{cases}
\sum_{r\in \hf+\Z_+}\psi^{+,i}_{-r}\psi^{-,j}_{r}, & \text{if $1\leq i, j\leq \ell$},\\
\sum_{r\in \hf+\Z_+}(-1)^{r+\hf}\psi^{+,i}_{-r}\psi^{+,j-\ell}_{r}, & \text{if $1\leq i\leq \ell$ and $\ell<j\leq n$},\\
\sum_{r\in \hf+\Z_+}(-1)^{r+\hf}\psi^{-,i-\ell}_{-r}\psi^{-,j}_{r}, & \text{if $\ell<i\leq n$ and $1\leq j\leq \ell$},\\
\sum_{r\in \hf+\Z_+}\psi^{-,i-\ell}_{-r}\psi^{+,j-\ell}_{r}, & \text{if $\ell<i,j\leq n$}.
\end{cases}
\end{equation}

For $n=2\ell+1$ and ${\rm G}_n={\rm O}_{2\ell+1}$, put $e_{ij}$ to be the same as in \eqref{eq:gl-1} except we assume that $\psi_{r}^{\pm,n}=\phi_r^{\pm}$.
Then one can check that \eqref{eq:gl-1} and \eqref{eq:gl-2} define an action of ${\rm GL}_n$ on $\mc{F}^{\frac{n}{2}}$ such that its restriction to ${\rm G}_n$ coincides with the action of ${\rm G}_n$ in Proposition \ref{duality}.
Recall that there exists an action of $\mf{l}_\infty\times {\rm GL}_{n}$ on $\mc{F}^{\frac{n}{2}}$ such that
\begin{equation}\label{eq:type A decomposition}
\mc{F}^{\frac{n}{2}}\cong \bigoplus_{\la\in \cP}L({\mf l}_\infty,\la')\otimes V^\la_{{\rm GL}_n}.
\end{equation}
On the other hand, by Proposition \ref{duality} we have
\begin{equation*}
\begin{split}
\mc{F}^{\frac{n}{2}}&\cong
\bigoplus_{\mu\in {\mc P}({\rm G}_n)}
L(\mf{g}_\infty,\La^{\mf g}(\mu))\otimes V_{{\rm G}_n}^\mu\\
&\cong \bigoplus_{\mu\in {\mc P}({\rm G}_n)}
\left(\bigoplus_{\la\in \cP} L(\mf{l}_\infty,\la')^{\oplus c^\mu_\la({\mf g})}\right)\otimes V_{{\rm G}_n}^\mu\\
&\cong \bigoplus_{\la\in \cP}  L(\mf{l}_\infty,\la')\otimes
\left(\bigoplus_{\mu\in {\mc P}({\rm G}_n)}(V_{{\rm G}_n}^\mu)^{\oplus c^\mu_\la({\mf g})}\right).
\end{split}
\end{equation*}
Combining with \eqref{eq:type A decomposition}, we have
\begin{equation*}
\left[ V^\la_{{\rm GL}_n} : V^\mu_{{\rm G}_n} \right]=c^\mu_\la({\mf g}).
\end{equation*}
The second equality follows from Corollary \ref{cor:stable branching}.
\qed

\begin{ex}\label{ex:GL to O-1}
{\rm Let $\la=(2,2)\in\cP$. The decomposition number $[V^\la_{{\rm GL}_n}: V^\mu_{{\rm O}_n}]$ for each $\mu\in \mc{P}({\rm O}_n)$ is equal to the number of ${\bf T}\in {\bf T}^{\mf d}(\mu,n)$ such that ${\bf T}\equiv_A H_{\la'}$. Using this fact, it is straightforward to check that
\begin{equation*}
[V^\la_{{\rm GL}_n}: V^\mu_{{\rm O}_n}]\neq 0 \ \Longleftrightarrow \
\mu\in 
\begin{cases}
\{\,(2,2,0,\ldots,0), (2,0,\ldots,0), (0,\ldots,0)\,\}, & \text{if $n\geq 4$},\\
\{\,(2,0,0), (0,0,0)\,\}, & \text{if $n=3$},\\
\{\,(0,0)\,\}, & \text{if $n=2$}.
\end{cases}
\end{equation*}

For $n=4$, we have ${\bf T}^{\mf d}(\mu,4)\subset \widehat{{\bf T}}^{\mf d}(\mu,4)={\bf T}^{\mf d}(\mu_1)\times {\bf T}^{\mf d}(\mu_2)$ for $\mu\in \mc{P}({\rm O}_4)$ with $[V^\la_{{\rm GL}_4}: V^\mu_{{\rm O}_4}]\neq 0$,
and the tableaux ${\bf T}=(T_1,T_2)\in {\bf T}^{\mf d}(\mu,4)$ such that ${\bf T}\equiv_A H_{\la'}$ are as follows: 
$$\resizebox{.9\hsize}{!}{$
\left(\
{\def\lr#1{\multicolumn{1}{|@{\hspace{.75ex}}c@{\hspace{.75ex}}|}{\raisebox{-.04ex}{$#1$}}}\raisebox{-.6ex}
{$\begin{array}{cc}
\\ \\
\cline{1-1}\cdashline{2-2}[0.5pt/1pt]
\lr{1}& \\ 
\cline{1-1} 
\lr{2}&\\ 
\cline{1-1}  
\end{array}$}}\quad , \quad
{\def\lr#1{\multicolumn{1}{|@{\hspace{.75ex}}c@{\hspace{.75ex}}|}{\raisebox{-.04ex}{$#1$}}}\raisebox{-.6ex}
{$\begin{array}{cc}
\\ \\
\cline{1-1}\cdashline{2-2}[0.5pt/1pt]
\lr{1}& \\ 
\cline{1-1} 
\lr{2}&\\ 
\cline{1-1}  
\end{array}$}}\
\right)\ \ \in {\bf T}^{\mf d}((2,2,0,0),4) \ \ \ \  
\left(\
{\def\lr#1{\multicolumn{1}{|@{\hspace{.75ex}}c@{\hspace{.75ex}}|}{\raisebox{-.04ex}{$#1$}}}\raisebox{-.6ex}
{$\begin{array}{cc}
\\ \\
\cline{1-1}\cdashline{2-2}[0.5pt/1pt]
\lr{1}& \\ 
\cline{1-1} 
\lr{2}&\\ 
\cline{1-1}  
\end{array}$}} \quad , \quad
{\def\lr#1{\multicolumn{1}{|@{\hspace{.75ex}}c@{\hspace{.75ex}}|}{\raisebox{-.04ex}{$#1$}}}\raisebox{-.6ex}
{$\begin{array}{cc}
\cline{2-2} 
&\lr{1} \\ 
\cline{2-2} 
&\lr{2} \\ 
\cline{2-2}\cdashline{1-2}[0.5pt/1pt]
\\
\\
\end{array}$}}\
\right)\ \ \in {\bf T}^{\mf d}((2,0,0,0),4) \ \ \ \
$}\\
$$\vskip 3mm

\ \ $\resizebox{.44\hsize}{!}{$
\left(\
{\def\lr#1{\multicolumn{1}{|@{\hspace{.75ex}}c@{\hspace{.75ex}}|}{\raisebox{-.04ex}{$#1$}}}\raisebox{-.6ex}
{$\begin{array}{cc}
\\ \\
\cdashline{1-2}[0.5pt/1pt]
 & \\ 
 & \\ 
\end{array}$}} \quad , \quad
{\def\lr#1{\multicolumn{1}{|@{\hspace{.75ex}}c@{\hspace{.75ex}}|}{\raisebox{-.04ex}{$#1$}}}\raisebox{-.6ex}
{$\begin{array}{cc}
\cline{1-2} 
\lr{1}&\lr{1} \\ 
\cline{1-2} 
\lr{2}&\lr{2} \\ 
\cline{1-2} 
\\
\\
\end{array}$}}\
\right)\ \ \in {\bf T}^{\mf d}((0,0,0,0),4) \ \ \ \
$}
$\vskip 3mm
\noindent where the dashed line denotes the line separating the body and the tail of $T_i$'s for $i=1,2$. Hence we have
\begin{equation*}
V^{(2,2)}_{{\rm GL}_4} = V^{(2,2,0,0)}_{{\rm O}_4} \oplus V^{(2,0,0,0)}_{{\rm O}_4} \oplus V^{(0,0,0,0)}_{{\rm O}_4}.
\end{equation*}
It is not difficult to see that the same decomposition of $V^{(2,2)}_{{\rm GL}_n}$ as above holds for $n\geq 4$, which is equal to the stable formula in Theorem \ref{thm:stable branching-dual}.

Let us consider the decomposition of $V^{\la}_{{\rm GL}_n}$ ($n=2,3$), where $\la$ is outside the stable range.
For $n=3$, we have ${\bf T}^{\mf d}(\mu,3)\subset \widehat{{\bf T}}^{\mf d}(\mu,3)={\bf T}^{\mf d}(\mu_1)\times {\bf T}^{\rm sp\,+}$ for $\mu\in \mc{P}({\rm O}_3)$ with $[V^\la_{{\rm GL}_3}: V^\mu_{{\rm O}_3}]\neq 0$,
and the tableaux ${\bf T}=(T_1,T_2)\in {\bf T}^{\mf d}(\mu,3)$ such that ${\bf T}\equiv_A H_{\la'}$ are 
$$\resizebox{.78\hsize}{!}{$
\left(\
{\def\lr#1{\multicolumn{1}{|@{\hspace{.75ex}}c@{\hspace{.75ex}}|}{\raisebox{-.04ex}{$#1$}}}\raisebox{-.6ex}
{$\begin{array}{cc}
\\ \\
\cline{1-1}\cdashline{2-2}[0.5pt/1pt]
\lr{1}& \\ 
\cline{1-1} 
\lr{2}&\\ 
\cline{1-1}  
\end{array}$}} \quad , \quad
{\def\lr#1{\multicolumn{1}{|@{\hspace{.75ex}}c@{\hspace{.75ex}}|}{\raisebox{-.04ex}{$#1$}}}\raisebox{-.6ex}
{$\begin{array}{c}
\cline{1-1} 
\lr{1} \\ 
\cline{1-1} 
\lr{2} \\ 
\cline{1-1}
\\
\\
\end{array}$}}\
\right)\ \ \in {\bf T}^{\mf d}((2,0,0),3) \ \ \ \
\left(\
{\def\lr#1{\multicolumn{1}{|@{\hspace{.75ex}}c@{\hspace{.75ex}}|}{\raisebox{-.04ex}{$#1$}}}\raisebox{-.6ex}
{$\begin{array}{cc}
\cline{2-2} 
&\lr{1} \\ 
\cline{2-2} 
&\lr{2} \\ 
\cline{2-2}\cdashline{1-1}[0.5pt/1pt]
 \\ \\ 
\end{array}$}} \quad , \quad
{\def\lr#1{\multicolumn{1}{|@{\hspace{.75ex}}c@{\hspace{.75ex}}|}{\raisebox{-.04ex}{$#1$}}}\raisebox{-.6ex}
{$\begin{array}{c}
\cline{1-1} 
\lr{1} \\ 
\cline{1-1} 
\lr{2} \\ 
\cline{1-1}
\\
\\
\end{array}$}}\
\right)\ \ \in {\bf T}^{\mf d}((0,0,0),3)
$}
$$\vskip 3mm

\noindent This implies that
\begin{equation*}
V^{(2,2)}_{{\rm GL}_3} = V^{(2,0,0)}_{{\rm O}_3} \oplus V^{(0,0,0)}_{{\rm O}_3}.
\end{equation*}
For $n=2$,  
$$\resizebox{.24\hsize}{!}{$
{\def\lr#1{\multicolumn{1}{|@{\hspace{.75ex}}c@{\hspace{.75ex}}|}{\raisebox{-.04ex}{$#1$}}}\raisebox{-.6ex}
{$\begin{array}{cc}
\cline{1-2} 
\lr{1}&\lr{1} \\ 
\cline{1-2} 
\lr{2}&\lr{2} \\ 
\cline{1-2} 
\end{array}$}}\
\ \in {\bf T}^{\mf d}((0,0),2) \ \ \ \
$}
$$
is the unique tableau in ${\bf T}^{\mf d}((0,0),2)$ such that ${\bf T}\equiv_A H_{\la'}$, and hence 
\begin{equation*}
V^{(2,2)}_{{\rm GL}_2} = V^{(0,0)}_{{\rm O}_2}.
\end{equation*}

}
\end{ex}

\begin{ex}{\rm Let $\la=(3,1)\in\cP$. As in Example \ref{ex:GL to O-1}, we can check that
\begin{equation*}
[V^\la_{{\rm GL}_n}: V^\mu_{{\rm O}_n}]\neq 0 \ \Longleftrightarrow \
\mu\in 
\begin{cases}
\{\,(3,1,0,\ldots,0), (2,0,\ldots,0), (1,1,0,\ldots,0)\,\}, & \text{if $n\geq 3$},\\
\{\,(2,0), (1,1)\,\}, & \text{if $n=2$}.
\end{cases}
\end{equation*}

For $n=4$, we have ${\bf T}^{\mf d}(\mu,4)\subset \widehat{{\bf T}}^{\mf d}(\mu,4)={\bf T}^{\mf d}(\mu_1)\times {\bf T}^{\mf d}(\mu_2)$ for $\mu\in \mc{P}({\rm O}_4)$ with $[V^\la_{{\rm GL}_4}: V^\mu_{{\rm O}_4}]\neq 0$,
and the tableaux ${\bf T}=(T_1,T_2)\in {\bf T}^{\mf d}(\mu,4)$ such that ${\bf T}\equiv_A H_{\la'}$ are as follows:
$$\resizebox{.9\hsize}{!}{$
\left(\
{\def\lr#1{\multicolumn{1}{|@{\hspace{.75ex}}c@{\hspace{.75ex}}|}{\raisebox{-.04ex}{$#1$}}}\raisebox{-.6ex}
{$\begin{array}{cc}
\\ 
\cline{1-1}\cdashline{2-2}[0.5pt/1pt]
\lr{1}& \\ 
\cline{1-1} 
\lr{2}&\\ 
\cline{1-1}  
\lr{3}&\\ 
\cline{1-1}  
\end{array}$}}\quad , \quad
{\def\lr#1{\multicolumn{1}{|@{\hspace{.75ex}}c@{\hspace{.75ex}}|}{\raisebox{-.04ex}{$#1$}}}\raisebox{-.6ex}
{$\begin{array}{cc}
\\ 
\cline{1-1}\cdashline{2-2}[0.5pt/1pt]
\lr{1}& \\ 
\cline{1-1} \\ \\
\end{array}$}}\
\right)\ \ \in {\bf T}^{\mf d}((3,1,0,0),4) \ \ \ \  
\left(\
{\def\lr#1{\multicolumn{1}{|@{\hspace{.75ex}}c@{\hspace{.75ex}}|}{\raisebox{-.04ex}{$#1$}}}\raisebox{-.6ex}
{$\begin{array}{cc}
\\ \\
\cline{1-1}\cdashline{2-2}[0.5pt/1pt]
\lr{1}& \\ 
\cline{1-1} 
\lr{3}&\\ 
\cline{1-1}  
\end{array}$}} \quad , \quad
{\def\lr#1{\multicolumn{1}{|@{\hspace{.75ex}}c@{\hspace{.75ex}}|}{\raisebox{-.04ex}{$#1$}}}\raisebox{-.6ex}
{$\begin{array}{cc}
\cline{2-2} 
&\lr{1} \\ 
\cline{2-2} 
&\lr{2} \\ 
\cline{2-2}\cdashline{1-2}[0.5pt/1pt]
\\
\\
\end{array}$}}\
\right)\ \ \in {\bf T}^{\mf d}((2,0,0,0),4) \ \ \ \
$}\\
$$\vskip 3mm

\ \ $\resizebox{.43\hsize}{!}{$
\left(\
{\def\lr#1{\multicolumn{1}{|@{\hspace{.75ex}}c@{\hspace{.75ex}}|}{\raisebox{-.04ex}{$#1$}}}\raisebox{-.6ex}
{$\begin{array}{cc}
\\ \\
\cline{1-1}\cdashline{1-2}[0.5pt/1pt]
\lr{1} & \\ 
\cline{1-1}
 & \\ 
\end{array}$}} \quad , \quad
{\def\lr#1{\multicolumn{1}{|@{\hspace{.75ex}}c@{\hspace{.75ex}}|}{\raisebox{-.04ex}{$#1$}}}\raisebox{-.6ex}
{$\begin{array}{cc}
\cline{2-2} 
&\lr{1} \\ 
\cline{2-2} 
&\lr{2} \\ 
\cline{1-2} 
\lr{3}&\\
\cline{1-1}
\\
\end{array}$}}\
\right)\ \ \in {\bf T}^{\mf d}((1,1,0,0),4) \ \ \ \
$}
$\vskip 3mm
\noindent Hence we have
\begin{equation*}
V^{(3,1)}_{{\rm GL}_4} = V^{(3,1,0,0)}_{{\rm O}_4} \oplus V^{(2,0,0,0)}_{{\rm O}_4} \oplus V^{(1,1,0,0)}_{{\rm O}_4}.
\end{equation*}
The above decomposition also holds for $n\geq 4$.

For $n=3$, we have ${\bf T}^{\mf d}(\mu,3)\subset \widehat{{\bf T}}^{\mf d}(\mu,3)={\bf T}^{\mf d}(\mu_1)\times {\bf T}^{\rm sp\,\pm}$ for $\mu\in \mc{P}({\rm O}_3)$ with $[V^\la_{{\rm GL}_3}: V^\mu_{{\rm O}_3}]\neq 0$,
and the tableaux ${\bf T}=(T_1,T_2)\in {\bf T}^{\mf d}(\mu,3)$ such that ${\bf T}\equiv_A H_{\la'}$ are as follows:
$$\resizebox{.75\hsize}{!}{$
\left(\
{\def\lr#1{\multicolumn{1}{|@{\hspace{.75ex}}c@{\hspace{.75ex}}|}{\raisebox{-.04ex}{$#1$}}}\raisebox{-.6ex}
{$\begin{array}{cc}
\\ 
\cline{1-1}\cdashline{2-2}[0.5pt/1pt]
\lr{1}& \\ 
\cline{1-1} 
\lr{2}&\\ 
\cline{1-1}  
\lr{3}&\\ 
\cline{1-1}  
\end{array}$}}\quad , \quad
{\def\lr#1{\multicolumn{1}{|@{\hspace{.75ex}}c@{\hspace{.75ex}}|}{\raisebox{-.04ex}{$#1$}}}\raisebox{-.6ex}
{$\begin{array}{c}
\\ 
\cline{1-1} 
\lr{1}  \\ 
\cline{1-1} \\ \\
\end{array}$}}\
\right)\ \ \in {\bf T}^{\mf d}((3,1,0),3)\ \ \ \
\left(\
{\def\lr#1{\multicolumn{1}{|@{\hspace{.75ex}}c@{\hspace{.75ex}}|}{\raisebox{-.04ex}{$#1$}}}\raisebox{-.6ex}
{$\begin{array}{cc}
\\ \\
\cline{1-1}\cdashline{2-2}[0.5pt/1pt]
\lr{1}& \\ 
\cline{1-1} 
\lr{3}&\\ 
\cline{1-1}  
\end{array}$}} \quad , \quad
{\def\lr#1{\multicolumn{1}{|@{\hspace{.75ex}}c@{\hspace{.75ex}}|}{\raisebox{-.04ex}{$#1$}}}\raisebox{-.6ex}
{$\begin{array}{c}
\cline{1-1} 
\lr{1} \\ 
\cline{1-1} 
\lr{2} \\ 
\cline{1-1}
\\
\\
\end{array}$}}\
\right)\ \ \in {\bf T}^{\mf d}((2,0,0),3) 
$}\\
$$\vskip 3mm

\quad\quad\ $\resizebox{.4\hsize}{!}{$
\left(\
{\def\lr#1{\multicolumn{1}{|@{\hspace{.75ex}}c@{\hspace{.75ex}}|}{\raisebox{-.04ex}{$#1$}}}\raisebox{-.6ex}
{$\begin{array}{cc}
\\ \\
\cline{1-1}\cdashline{1-2}[0.5pt/1pt]
\lr{1} & \\ 
\cline{1-1}
 & \\ 
\end{array}$}} \quad , \quad
{\def\lr#1{\multicolumn{1}{|@{\hspace{.75ex}}c@{\hspace{.75ex}}|}{\raisebox{-.04ex}{$#1$}}}\raisebox{-.6ex}
{$\begin{array}{c}
\cline{1-1} 
\lr{1} \\ 
\cline{1-1} 
\lr{2} \\ 
\cline{1-1} 
\lr{3}\\
\cline{1-1}
\\
\end{array}$}}\
\right)\ \ \in {\bf T}^{\mf d}((1,1,0),3) \ \ \ \
$}
$\vskip 3mm

\noindent Hence we have
\begin{equation*}
V^{(3,1)}_{{\rm GL}_3} = V^{(3,1,0)}_{{\rm O}_3} \oplus V^{(2,0,0)}_{{\rm O}_3} \oplus V^{(1,1,0)}_{{\rm O}_3}.
\end{equation*}
Finally, for $n=2$, the tableaux $T\in {\bf T}^{\mf d}(\mu,2)$ such that $T\equiv_A H_{\la'}$ are given by    

$$\resizebox{.53\hsize}{!}{$
{\def\lr#1{\multicolumn{1}{|@{\hspace{.75ex}}c@{\hspace{.75ex}}|}{\raisebox{-.04ex}{$#1$}}}\raisebox{-.6ex}
{$\begin{array}{cc}
\cline{2-2} 
&\lr{1} \\ 
\cline{2-2} 
&\lr{2} \\ 
\cline{1-2}
\lr{1}& \\
\cline{1-1}
\lr{3}\\
\cline{1-1}
\end{array}$}}\ \ \in {\bf T}^{\mf d}((2,0),2) \ \ \ \ \ \ \ \ \ \
{\def\lr#1{\multicolumn{1}{|@{\hspace{.75ex}}c@{\hspace{.75ex}}|}{\raisebox{-.04ex}{$#1$}}}\raisebox{-.6ex}
{$\begin{array}{cc}
\cline{2-2} 
&\lr{1} \\ 
\cline{2-2} 
&\lr{2} \\ 
\cline{1-2} 
\lr{1}&\lr{3}\\
\cline{1-2}
\\
\end{array}$}}\ \ \in {\bf T}^{\mf d}((1,1),2) \ \ \ \ \ \ \ \
$}\\
$$\vskip 3mm
\noindent and hence
\begin{equation*}
V^{(3,1)}_{{\rm GL}_2} = V^{(2,0)}_{{\rm O}_2} \oplus V^{(1,1)}_{{\rm O}_2}.
\end{equation*}

}
\end{ex}

\begin{rem}{\rm
Let $\la^{(1)},\ldots,\la^{(r)}\in \cP$ with $\ell(\la^{(i)})\leq n$ for $1\leq i\leq r$. For $\la\in \cP$ with $\ell(\la)\leq n$, let 
\begin{equation*}
c^{\la}_{\la^{(1)}\cdots \la^{(r)}}=
\left[ 
V_{{\rm GL}_n}^{\la^{(1)}}\otimes \cdots \otimes V_{{\rm GL}_n}^{\la^{(r)}} : V_{{\rm GL}_n}^\la
\right].
\end{equation*}
Then for $\mu\in {\mc P}({\rm G}_n)$, we have 
\begin{equation}\label{eq:branching multiple tensor}
\left[ 
V_{{\rm GL}_n}^{\la^{(1)}}\otimes \cdots \otimes V_{{\rm GL}_n}^{\la^{(r)}} : V^{\mu}_{{\rm G}_n}
\right]
=\sum_{\la\in \cP} \left[ V^\la_{{\rm GL}_n} : V^\mu_{{\rm G}_n} \right]c^{\la}_{\la^{(1)}\cdots \la^{(r)}}.
\end{equation}
In particular, if $|\la^{(1)}|+\cdots+|\la^{(r)}|\leq \frac{n}{2}$, where $|\la^{(i)}|$ denotes the sum of non-zero parts of $\la^{(i)}$, then $\ell(\la)\leq \frac{n}{2}$ and hence \eqref{eq:branching multiple tensor} becomes
\begin{equation}\label{eq:branching multiple tensor-2}
\left[ 
V_{{\rm GL}_n}^{\la^{(1)}}\otimes \cdots \otimes V_{{\rm GL}_n}^{\la^{(r)}} : V^{\mu}_{{\rm G}_n}
\right]
=\sum_{\la\in \cP}\sum_{\delta\in \cP_{\mf g}} c^{\la}_{\delta'\mu}c^{\la}_{\la^{(1)}\cdots \la^{(r)}}.
\end{equation}
By using Theorem \ref{thm:stable branching-dual}, we obtain a new combinatorial formula for \eqref{eq:branching multiple tensor} or a combinatorial extension of \eqref{eq:branching multiple tensor-2} to arbitrary $\la^{(i)}$'s. We should remark that a $q$-analogue of \eqref{eq:branching multiple tensor} is introduced in \cite{Le07}, and its stable limit is closely related with the energy function on a tensor product of finite affine crystals (see for example \cite{LOS}). It would be interesting to find a combinatorial interpretation of the $q$-analogue of \eqref{eq:branching multiple tensor} or \eqref{eq:branching multiple tensor-2} in terms of spinor model.
}
\end{rem}

\section{Character of holomorphic discrete series}\label{sec:oscillator}

\subsection{Superization of ${\bf T}^{\mf g}(\la,n)$}
Let $\mc{A}$ be a linearly ordered countable set with a $\mathbb{Z}_2$-grading $\mc{A}=\mc{A}_0\sqcup\mc{A}_1$.
For a skew Young diagram $\lambda/\mu$, we define the notion of $\mc{A}$-semistandard tableaux of shape $\lambda/\mu$ as usual,  that is, (1) the entries in each row (resp. column) are
weakly increasing from left to right (resp. from top to bottom), (2)
the entries in $\mc{A}_0$ (resp. $\mc{A}_1$) are strictly increasing in each
column (resp. row).
We denote by $SST_{\mc{A}}(\lambda/\mu)$ the set of $\mc{A}$-semistandard tableaux of shape $\lambda/\mu$. Let $x_{\mc{A}}=\{\,x_a\,|\,a\in\mc{A}\,\}$ be the set of formal commuting variables indexed by $\mc{A}$. For $\lambda\in\cP$, let $s_{\lambda}(x_{\mc{A}})=\sum_{T}x_\mc{A}^T$ be the super Schur function corresponding to $\lambda$, where the sum is over $T\in SST_{\mc{A}}(\lambda)$ and $x_\mc{A}^T=\prod_{a}x_a^{m_a}$ with $m_a$ the number of occurrences of $a$ in $T$.

For $\la\in {\mc P}({\rm G}_n)$, one can define ${\bf T}^{\mf g}_{\mc A}(\la,n)$ by replacing ($\N$-)semistandard tableaux with $\mc A$-semistandard tableaux in Definition \ref{def:psst}. (See \cite[Section 6]{K13} and \cite[Section 3]{K14} for more details, where ${\bf T}^{\mf g}_{\mc A}(\la,n)$ is denoted by ${\bf T}^{\mf g}_{\mc A}(\la,n/\epsilon)$.)
When $\mc{A}$ is (a subset of) $\N$, we assume that $\mc{A}_0=\mc{A}$ and it is equipped with the usual linear ordering on $\N$. In particular, we have ${\bf T}^{\mf g}_{\N}(\la,n)={\bf T}^{\mf g}(\la,n)$.

We define the character of ${\bf T}^{\mf g}_{\mc A}(\la,n)$ by
\begin{equation*}
{\rm ch}{\bf T}^{\mf g}_{\mc A}(\la,n) =z^{n/\epsilon}\sum_{{\bf T}}\prod_{i\geq 1} x_\mc{A}^{T_i},
\end{equation*}
where $z$ is another formal variable and the sum is over ${\bf T}=(T_1,T_2,\ldots)\in {\bf T}^{\mf g}_{\mc A}(\la,n)$. 
It is shown that under certain choices of $\mc A$ (more precisely, for $\mc{A}$ with $0<|\mc{A}_0|<\infty$), $\{\,{\rm ch}{\bf T}^{\mf g}_{\mc A}(\la,n)\,|\,\la\in \bigsqcup_n\mc{P}({\rm G}_n)\,\}$ gives the characters of a family of irreducible representation of ortho-symplectic Lie superalgebra in a semisimple tensor category \cite[Theorem 7.6]{K13} and \cite[Theorem 4.6]{K14}.

\subsection{Unitarizable representations}

From now on, we assume that $2\leq k\leq \infty$ and $({\mf{g}}_k,{\rm G}_n)$ denotes one of the following pairs
\begin{equation}\label{eq:dual pairs-2}
({\mf d}_k, {\rm Sp}_n),\ \ ({\mf b}^\bullet_k, {\rm Pin}_n),\ \ ({\mf b}^\bullet_k, {\rm Spin}_n),\ \  ({\mf c}_k, {\rm O}_n).
\end{equation} 
Put
\begin{equation*}
{\mc P}({\rm G}_n)_k=\{\,\la\,|\,\la\in {\mc P}({\rm G}_n),\ \ell(\la)\leq k\,\}.
\end{equation*}
For $\la\in {\mc P}({\rm G}_n)_k$, we define a weight $\La^{\mf g}_-(\la)$ for ${\mf{g}}_k$ by
\begin{equation*}\label{eq:hw-2}
\Lambda^{\mf g}_-(\la)=
-\frac{n}{\epsilon}\La^{\mf g}_0+\la_1\epsilon_1+\la_2\epsilon_2+\cdots.
\end{equation*}
Note that $\La^{\mf g}_0=-(\epsilon/2)(\epsilon_1+\cdots+\epsilon_k)$ for $k<\infty$.
We denote by $L({\mf{g}}_k, \La^{\mf g}_-(\la))$ the irreducible highest weight ${\mf{g}}_k$-module with highest weight $\La^{\mf g}_-(\la)$.
Let $\mf h_k$ be the Cartan subalgebra of $\g_k$. We define the character of $L({\mf{g}}_k, \La^{\mf g}_-(\la))$ by ${\rm ch}L({\mf{g}}_k, \La^{\mf g}_-(\la))=\sum_{\mu}\dim L({\mf{g}}_k, \La^{\mf g}_-(\la))_\mu e^\mu$, where $L({\mf{g}}_k, \La^{\mf g}_-(\la))_\mu$ is the $\mu$- weight space and $e^{\mu}$ is the basis element of the group algebra $\mathbb{Q}[{\mf h}^*_k]$. 

It is well-known that $L({\mf{g}}_k, \La^{\mf g}_-(\la))$ for $\la\in\mc{P}({\rm G}_n)_k$ form a family of infinite-dimensional unitarizable representations, which appear in the decomposition of the Segal-Shale representation of the metaplectic group with respect to the reductive dual pair $(\mf{g}_k,{\rm G}_n)$ (see \cite{HTW} and many related references therein for $\mf{g}=\mf{c}, \mf{d}$ and $k<\infty$, and also \cite{CKW,Wa} for $\mf{g}=\mf{b}^\bullet$ or $k=\infty$). For $\la\in {\mc P}({\rm G}_n)_k$, $L({\mf{g}}_k, \La^{\mf g}_-(\la))$ is not an integrable $\mf{g}_k$-module but it decomposes into a sum of integrable highest weight $\mf{l}_k$-modules. For $\mu\in {\mc P}({\rm G}_m)_k$ and $\nu\in {\mc P}({\rm G}_n)_k$, the tensor product $L({\mf{g}}_k, \La^{\mf g}_-(\mu))\otimes L({\mf{g}}_k, \La^{\mf g}_-(\nu))$ is completely reducible and is a direct sum of $L({\mf{g}}_k, \La^{\mf g}_-(\la))$'s for $\la\in {\mc P}({\rm G}_{m+n})_k$.

Let us assume that $(\mf{g}^\vee_k,{\mf g}_k)$ denotes one of the pairs
\begin{equation}\label{eq:pair of Lie}
(\mf{c}_k,\mf{d}_k), \ (\mf{b}_k,\mf{b}^\bullet_k), \ (\mf{d}_k, \mf{c}_k)
\end{equation}
such that $(\mf{g}^\vee_\infty,{\rm G}_n)$ and $({\mf g}_\infty,{\rm G}_n)$ are the dual pairs with the same ${\rm G}_n$ in \eqref{eq:dual pairs} and \eqref{eq:dual pairs-2}, respectively. It was recently observed from the theory of super duality \cite{CLW} that there is an equivalence between two parabolic BGG categories of modules over $\mf{g}^\vee_\infty$ and ${\mf g}_\infty$ with respect to $\mf{l}_\infty$, which sends $L({\mf{g}}^\vee_\infty,\La^{\mf g^\vee}(\la))$ to $L({\mf{g}}_\infty,\La^{{\mf g}}_-(\la))$ for $\la\in \bigsqcup_n\mc{P}({\rm G}_n)$.

Let
$\N'=\{\,p'\,|\,p\in\N \,\}$ be the set of primed positive integers of odd degree with a linear ordering $<$ where $p'<q'$ if and only if $p<q$.
For $\la\in {\mc P}({\rm G}_n)_k$, we put
\begin{equation*}\label{eq:tableaux_oscillator}
\begin{split}
{\mathbb T}^{{\mf{g}}}_k(\la,n)=
\begin{cases}
{\bf T}^{{\mf{g}}^\vee}_{\N'}(\la,n),& \text{if $k=\infty$},\\
{\bf T}^{{\mf{g}}^\vee}_{[k]'}(\la,n),& \text{if $k<\infty$},
\end{cases}
\end{split}
\end{equation*}
where $[k]'=\{\,1',2',\ldots,k'\,\}\subset \N'$.
Now, we are in a position to give a combinatorial character formula for $L({\mf{g}}_k, \La^{\mf g}_-(\la))$.

\begin{thm}\label{thm:char of unitary}
For $\la\in \mc{P}({\rm G}_n)_k$, we have
\begin{equation*}
{\rm ch}L({\mf{g}}_k,\La^{\mf g}_-(\la)) = {\rm ch}{\mathbb T}^{{\mf{g}}}_k(\la,n),
\end{equation*}
where $z=e^{-\La_0^{\mf g}}$ and $x_{i'}=e^{\epsilon_i}$ for $i\geq 1$.
\end{thm}
\pf The proof is almost the same as in \cite[Theorem 7.6]{K13} and \cite[Theorem 4.6]{K14}.
First consider the case when $k=\infty$. Let $(\mf{g}^\vee_k,{\mf g}_k)$ be as in \eqref{eq:pair of Lie}.
It follows from \cite[Theorem 4.6]{CLW} (see also \cite[Section 2.3]{CLW}) that if
\begin{equation}\label{eq:Schur comb}
z^{n/\epsilon}{\rm ch}L({\mf{g}}^\vee_\infty,\La^{\mf g^\vee}(\la))=\sum_{\mu\in \cP}K_{\la\mu}s_{\mu}(x_{\N})
\end{equation}
for some non-negative integers $K_{\la\mu}$, then
\begin{equation}\label{eq:Schur comb-2}
z^{-n/\epsilon}{\rm ch}L({\mf{g}}_\infty,\La^{{\mf g}}_-(\la))
=\sum_{\mu\in \cP}K_{\la\mu}s_{\mu'}(x_{\N})=\sum_{\mu\in \cP}K_{\la\mu}s_{\mu}(x_{\N'}).
\end{equation}
Hence by \cite[Theorem 6.12, Corollary 6.13]{K13} and \cite[Theorem 3.9]{K14}, we have
\begin{equation*}\label{eq:Schur comb-3}
{\rm ch}L({\mf{g}}_\infty,\La^{{\mf g}}_-(\la)) ={\rm ch}{\bf T}^{{\mf{g}}^\vee}_{\N'}(\la,n)= {\rm ch}{\mathbb T}^{{\mf{g}}}_\infty(\la,n).
\end{equation*}
Since ${\rm ch}L({\mf{g}}_k,\La^{{\mf g}}_-(\la))$ for $k<\infty$ is obtained by specializing $x_{i'}=0$ for $i>k$ in ${\rm ch}L({\mf{g}}_\infty,\La^{{\mf g}}_-(\la))$ \cite[Lemma 3.2]{CLW}, we have 
${\rm ch}L({\mf{g}}_k,\La^{\mf g}_-(\la)) = {\rm ch}{\mathbb T}^{{\mf{g}}}_k(\la,n)$.
\qed\vskip 2mm

As an immediate corollary, we have the following tensor product and branching multiplicity formula.

\begin{cor}\label{cor:branching for unitary}
\mbox{}
\begin{itemize}
\item[(1)] For $\mu\in {\mc P}({\rm G}_m)_k$, $\nu\in {\mc P}({\rm G}_n)_k$, and $\la\in \mc{P}({\rm G}_{m+n})_k$, we have
\begin{equation*}
c^{\la}_{\mu\nu }({\mf g}^\vee)=
\left[
L({\mf g}_k,\Lambda^{\mf g}_-(\mu))\otimes L({\mf g}_k,\Lambda^{\mf g}_-(\nu)):L({\mf g}_k,\Lambda^{\mf g}_-(\la))
\right].
\end{equation*}

\item[(2)] For $\mu\in \mc{P}({\rm G}_n)_k$ and $\la\in \cP$ with $\ell(\la)\leq k$,
\begin{equation*}
c_{\la}^\mu(\mf{g}^\vee)=[L({\mf{g}}_k,\La^{\mf g}_-(\mu)): L({\mf l}_k,\la)].
\end{equation*}

\end{itemize}
\end{cor}

\pf Let $\omega$ denote the involution on the ring of symmetric functions sending $s_\la(x_{\N})$ to $s_{\la'}(x_{\N})$ for $\la\in \cP$.
We see from \eqref{eq:Schur comb} and \eqref{eq:Schur comb-2} that 
$z^{-n/\epsilon}{\rm ch}L({\mf{g}}_\infty,\La^{{\mf g}}_-(\la))=\omega\left(z^{n/\epsilon}{\rm ch}L({\mf{g}}^\vee_\infty,\La^{{\mf g}^\vee}(\la))\right)$. Hence (1) and (2) for $k=\infty$ follow from Theorem \ref{thm:char of unitary}, Propositions \ref{prop:tensor product} and \ref{prop:branching to levi}. 
The proof for the case when $k<\infty$ is obtained by specializing $x_{i'}=0$ for $i>k$ in the characters.
\qed

\subsection{Characters in a stable range}
Fix $k$. Let $\mf{g}^+_k$ be the Borel subalgebra of $\mf g_k$ with respect to the simple roots in Section \ref{subsec:notations}, and let $\mf p_k=\mf l_k+\mf g^+_k$ be the parabolic subalgebra associated to $\mf l_k$. For $\la\in \mc{P}({\rm G}_n)_k$, let $V(\mf g_k,\La^{\mf g}_-(\la))={\rm Ind}^{\mf{g}_k}_{\mf p_k}L(\mf{l}_k,\La^{\mf g}_-(\la))$ be the generalized Verma module with highest weight $\La^{\mf g}_-(\la)$, whose maximal irreducible quotient is $L({\mf{g}}_k, \La^{\mf g}_-(\la))$.

In this subsection, we prove that $L({\mf{g}}_k, \La^{\mf g}_-(\la))$ is equal to $V({\mf{g}}_k, \La^{\mf g}_-(\la))$, and hence its character has a nice product form when $\la$ is in a certain stable range (these modules are referred to as {\em  holomorphic discrete series}). It is an already well-known result when $\mf{g}=\mf{c}$, $\mf{d}$, and plays a crucial role in establishing various families of stable branching rules for classical groups in a unified way \cite{HTW}.
Nevertheless the proof here is given in a completely different and purely combinatorial way by decomposing ${\mathbb T}^{{\mf{g}}}_k(\la,n)$ into $\mf{l}_k$-crystals. Moreover the result for $\mf{g}=\mf{b}^\bullet$ obtained here seems to be new.

We put $[k]=\{\,1,\ldots,k\,\}\subset \N$ as a $\Z_2$-graded linearly ordered set, which can be viewed as the crystal of the natural representation of $\mf{l}_k$. As in Section \ref{subsec:tableaux}, one can regard $SST_{[k]}(\la/\mu)$ as an $\mf{l}_k$-crystal for a skew Young diagram $\la/\mu$. For $T\in SST_{\mc A}(\la/\mu)$, we define $T'$ to be the tableau of shape $\la'/\mu'$ obtained by flipping $T$ with respect to the main diagonal and replacing $i$ (resp. $i'$) with $i'$ (resp. $i$) when $\mc A$ is  $[k]$ (resp. $[k]'$).

Let $\la\in \mc{P}({\rm G}_n)_k$ be given. We regard ${\mathbb T}^{{\mf{g}}}_k(\la,n)$ as an $\mf{l}_k$-crystal by identifying ${\mathbb T}=(T_1,\ldots,T_r)\in {\mathbb T}^{{\mf{g}}}_k(\la,n)$ with $(T_1)'\otimes \cdots \otimes (T_r)'$.
For $\mathbb{T}=(T_1,\ldots,T_r)\in {\mathbb T}^{{\mf{g}}}_k(\la,n)$, 
we put $w({\mathbb T})= w(T_r)\ldots w(T_1)$, and define $L({\mathbb T})$ to be the maximal length of a strictly decreasing subword of $w({\mathbb T})$. It is not difficult to see that if ${\mathbb{T}}\equiv_A S$ for some $S\in SST_{[k]}(\mu)$, then $L(\mathbb{T})=\ell(\mu)$.

\begin{lem}\label{lem:aux-4}
We have $k\geq L({\mathbb T})\geq \ell(\la)$.
\end{lem}
\pf Note that ${\mathbb T}^{{\mf{g}}}_k(\la,n)$ can be defined in the same way as in ${\bf T}^{{\mf{g}}}(\la,n)$ except that we replace all $\leq$ with $<$ in Definition \ref{admissible} since the elements in $[k]'$ are of odd degree. Then
Lemma \ref{lem:aux-1} also holds for $(T,S)\in {\mathbb T}^{\mf g}_k(a)\times {\mathbb T}^{\mf g}_k(a')$ $(a\geq a')$ such that $T\prec S$, where $\leq$ is also replaced with $<$. This implies that there exists a strictly decreasing subword $T^{\tt tail}_{\ell(\la)}(1)\cdots T^{\tt tail}_{1}(1)$, and hence $L(\mathbb{T})\geq \ell(\la)$. On the other hand, we have $k\geq L({\mathbb T})$ since the entries of ${\mathbb T}$ are  from $[k]'$.
\qed \vskip 2mm

Let us define ${\mathbb T}^{\texttt{tail}}$ and ${\mathbb T}^{\texttt{body}}$ be defined in the same way as in \eqref{eq:body and tail}.

\begin{lem}\label{lem:aux-5}
If $L({\mathbb T})\leq \frac{n}{2}$, then we have
\begin{itemize}
\item[(1)] ${\mathbb T}^{\tt tail}\in SST_{[k]'}(\la')$, and ${\mathbb T}^{\tt body}\in SST_{[k]'}(\mu^\pi)$ for some $\mu\in\cP_{\mf g^\vee}$,

\item[(2)] ${\mathbb T}\equiv_A ({\mathbb T}^{\tt tail})'\otimes ({\mathbb T}^{\tt body})'$.
\end{itemize}
\end{lem}
\pf The proof is almost parallel to that of Lemma \ref{lem:aux-3}.
\qed
\vskip 2mm

Now, we have the following characterization of ${\mathbb T}^{{\mf{g}}}_k(\la,n)$ as an $\mf{l}_k$-crystal in a stable range, which is the main result in this section.

\begin{thm}\label{thm:stable oscillator}
Let $\la\in \mc{P}({\rm G}_n)_k$ be given. If $n\geq 2k$, then we have an isomorphism of $\mf{l}_k$-crystals
\begin{equation*}
\begin{split}
& {\mathbb T}^{{\mf{g}}}_k(\la,n) \longrightarrow \bigsqcup_{\substack{\mu\in \cP_{\mf g} \\ \ell(\mu)\leq k}} SST_{[k]}(\la) \otimes SST_{[k]} (\mu^\pi), \\
\end{split}
\end{equation*}
which maps ${\mathbb T}$ to $({\mathbb T}^{\tt tail})'\otimes ({\mathbb T}^{\tt body})'$.
\end{thm}
\pf 
Let $\mathbb{B}$ denote the right-hand side of the above isomorphism.
Since $2k\leq n$, we have
$L(\mathbb{T})\leq k \leq \frac{n}{2}$
for all ${\mathbb T}\in {\mathbb T}^{{\mf{g}}}_k(\la,n)$. By Lemma \ref{lem:aux-5}, the map ${\mathbb T}\mapsto ({\mathbb T}^{\tt tail})'\otimes ({\mathbb T}^{\tt body})'$ is a well-defined morphism of $\mf{l}_k$-crystals from ${\mathbb T}^{{\mf{g}}}_k(\la,n)$ to $\mathbb{B}$. The bijectiveness of the map can be proved by similar arguments as in Theorem \ref{thm:stable tensor}.
\qed\vskip 2mm

We recover the following stability of ${\rm ch}L({\mf{g}}_k, \La^{\mf g}_-(\la))$ in a purely combinatorial way (cf.~\cite[Theorem 3.2]{HTW} and references therein).

\begin{cor}\label{cor:stable branching for unitary}
Let $\la\in \mc{P}({\rm G}_n)_k$ be given. If $n\geq 2k$, then we have
\begin{equation*}
{\rm ch}L({\mf{g}}_k, \La^{\mf g}_-(\la))=
z^{n/\epsilon}\Delta^{\mf g}(x_{[k]}) s_\la(x_{[k]}),
\end{equation*}
where
\begin{equation*}
\Delta^{\mf g}(x_{[k]})=
\begin{cases}
\prod_{1\leq i\leq k}(1-x_i)^{-1}\prod_{1\leq i<j\leq k}(1-x_ix_j)^{-1}, & \text{if $\mf{g}=\mf{b}^\bullet$}, \\
\prod_{1\leq i\leq k}(1-x_i^2)^{-1}\prod_{1\leq i<j\leq k}(1-x_ix_j)^{-1}, & \text{if $\mf{g}=\mf{c}$}, \\
\prod_{1\leq i<j\leq k}(1-x_ix_j)^{-1}, & \text{if $\mf{g}=\mf{d}$}.
\end{cases}
\end{equation*}
In particular, we have $L({\mf{g}}_k, \La^{\mf g}_-(\la))=V({\mf{g}}_k, \La^{\mf g}_-(\la))$.
\end{cor}
\pf It follows from Theorem \ref{thm:stable oscillator} and the well-known identity (cf.~\cite{Mac})
\begin{equation*}
\Delta^{\mf g}(x_{[k]})=\sum_{\substack{\la\in \cP_{\mf g} \\ \ell(\mu)\leq k}}s_\mu(x_{[k]}) 
\end{equation*}
that ${\rm ch}L({\mf{g}}_k, \La^{\mf g}_-(\la))=
z^{n/\epsilon}\Delta^{\mf g}(x_{[k]}) s_\la(x_{[k]})$.
Let $\mf u^-_k$ be the opposite nilradical of $\mf p_k$, and $U(\mf u^-_k)$ its universal enveloping algebra. We have ${\rm ch}U(\mf u^-_k)=\Delta^{\mf g}(x_{[k]})$, and hence ${\rm ch}V({\mf{g}}_k, \La^{\mf g}_-(\la))=z^{n/\epsilon}\Delta^{\mf g}(x_{[k]}) s_\la(x_{[k]})$. 
Since $L({\mf{g}}_k, \La^{\mf g}_-(\la))$ is a maximal quotient of $V({\mf{g}}_k, \La^{\mf g}_-(\la))$, we have $L({\mf{g}}_k, \La^{\mf g}_-(\la))=V({\mf{g}}_k, \La^{\mf g}_-(\la))$.
\qed

\begin{rem}{\rm
One can deduce Corollary \ref{cor:stable branching for unitary} from Theorem \ref{thm:stable branching}, Corollary \ref{cor:branching for unitary} (2), and the Littlewood identity (cf.~\cite{Mac}). Recall that the stability of character ${\rm ch}L({\mf{g}}_k, \La^{\mf g}_-(\la))$ in Corollary \ref{cor:stable branching for unitary} plays a crucial role in obtaining the stable tensor product multiplicities for $L({\mf{g}}_k, \La^{\mf g}_-(\la))$ in \cite{HTW}. We can also recover this stable tensor product multiplicity directly from Corollaries \ref{cor:stable tensor} and \ref{cor:branching for unitary} (1).

}
\end{rem}

\end{document}